\documentclass{article}

\usepackage{amssymb}
\usepackage[latin9]{inputenc}

\providecommand{\keywords}[1]
{
  \small	
  \textbf{\textit{Keywords---}} #1
}

\title{From Tractatus to Later Writings and Back - \\
New Implications from the \emph{Nachlass}\footnote{Preprint. Final version to appear in \emph{SATS -- Northern European Journal of Philosophy,} Walter de Gruyter.}}

\author{Ruy J.G.B. de Queiroz\\
Centro de Inform\'atica\\
Universidade Federal de Pernambuco\\
Recife, Brazil\\
ruy@cin.ufpe.br}

\date{}

\begin{document}

\maketitle

\begin{abstract}
As a celebration of the \emph{Tractatus} 100th anniversary it might be worth revisiting its relation to the later writings. From the former to the latter, David Pears recalls that ``everyone is aware of the holistic character of Wittgenstein's later philosophy, but it is not so well known that it was already beginning to establish itself in the \emph{Tractatus}" (\emph{The False Prison}, 1987).
From the latter to the former, Stephen Hilmy's (\emph{The Later Wittgenstein}, 1987) extensive study of the \emph{Nachlass} has helped removing
classical misconceptions such as Hintikka's claim that ``Wittgenstein in the \emph{Philosophical Investigations} almost completely gave up the calculus analogy." Hilmy points out that even in the \emph{Investigations} one finds the use of the calculus/game paradigm to the understanding of language, such as ``in operating with the word" (Part I, \S 559) and ``it plays a different part in the calculus". Hilmy also quotes from a late (1946) unpublished manuscript (MS 130) ``this sentence has use in the calculus of language"), which seems to be compatible with ``asking whether and how a proposition can be verified is only a particular way of asking `How do you mean?'." Central in this back and forth there is an aspect which seems to deserve attention in the discussion of a semantics for the language of mathematics which might be based on (normalisation of) proofs and/or Hintikka/Lorenzen game-dialogue: the explication of consequences. Such a discussion is substantially supported by the use of the open and searchable \emph{The Wittgenstein Archives at the University of Bergen}. These findings are framed within the discussion of the meaning of logical constants in the context of natural deduction style rules of inference.\footnote{The material reported herein (partly presented at the online meeting, 14--17 Sept 2021, \emph{100 anos do Tractatus Logico-Philosophicus}) further develops an ongoing project of putting forward a proposal of giving a formal counterpart to the `meaning is use' dictum advocated intensively by the later Wittgenstein, at the same time attempting to find a common thread in Wittgenstein early as well as later writings. For that reason it is based on a series of previous publications from which some parts are reused here to reinforce the flow of the general argument: \emph{Wittgenstein Symposium} (1988b,1989,1992), \emph{Dialectica} (1988a,1991,1994), \emph{Logic J.\ of the IGPL} (2001), \emph{Studia Logica} (2008).  
We believe that there is novelty in the present paper not least because we have made extensive use of Wittgenstein's \emph{Nachlass} which was not done in any of the previous papers. Thus, as a new step in a series of essays going back a few decades, it certainly contains some overlap with the previous publication intended to bring context and self-containedness to the present manuscript, this substantially reinforces earlier arguments. (Gratitude to the editors of \emph{Wittgenstein Archives at the University of Bergen} (\emph{WAB}) who provided Wittgenstein scholarship with such a valuable gift as the \emph{Nachlass} in free online and searchable form!)}
\end{abstract}


\section{Introduction}
In order to make sense of what we claim constitutes a common thread of Wittgenstein's view on the connections between meaning, use and consequences, going from the \emph{Tractatus} to later writings and back, and take this as the basis for a proposal for a formal counterpart of a `meaning-as-use' (dialogical/game-theoretical) semantics of the language of predicate logic, we shall need to be concise, and yet to bring in the key excerpts from Wittgenstein oeuvre (including the \emph{Nachlass})  and from those formal semanticists who defend a different interpretation of the connections between proofs and meaning. The aim is to consider the so-called rules of proof reduction as a formal counterpart to the explanation of the (immediate) consequences of a proposition. This contrasts markedly to a different perspective from the verificationist theories of meaning as put forward by Gentzen, Dummett, Prawitz, Martin-L\"of and many others, suggesting an approach which has a `pragmatist' slant to the semantics of predicate logic. Accordingly, we consider several passages from Wittgenstein's published as well as unpublished writings to build a whole picture of a formal counterpart to `meaning is use'.

The legendary contrast between `early' and `late' Wittgenstein has been countered by such scholars as Pears (1987) and Hilmy (1987). One important, yet neglected continuity in Wittgenstein's central issues and proposals concerns `meaning' as `use'. Linking his early and later texts on `meaning as use' is his appeal to direct consequences of a term or phrase, reflected e.g.\ in his speaking of language as a `calculus'. These passages are crucial to Wittgenstein's view of `meaning as use', though they have been widely neglected in scholarly literature. The centrality and significance of these passages are corroborated and augmented here by renewed examination of Wittgenstein's \emph{Nachla\ss}, as a further step in the present author's continuing research on this topic. Explicating these passages reveals close and important parallels to recent dialogue semantics, game semantics, and to Peirce's attention to the complementary perspectives of speakers and interpreters. Reduction rules, introduction rules, elimination rules and game-theoretic accounts of logical connectives form sets of complementary functional specifications, which are central to Wittgenstein's own account of `meaning as use'. These issues are introduced in section 1, in which Wittgenstein's abiding interest in the direct consequences of terms or phrases are documented and examined. Section 2 documents and examines the significance of their pervasive disregard in the literature. Additionally, section 2  examines Wittgenstein's continuing, central interest in the direct consequences of terms or phrases in later publications. In this connection, section 3 examines two current approaches to the interface of meaning, knowledge and logic: Lorenzen's dialogue semantics and Hintikka's game semantics, and shows how these converge in a common set of inference schemata.

The intuition regarding the connections between meaning and use (as usefulness), which appears explicitly in Wittgenstein's very late remarks such as ``Meaning, function, purpose, usefulness -- interconnected concepts"\footnote{\S 291, p.\ 41$^\mathrm{e}$, {\em Last Writings on the Philosophy of Psychology}}, and underlies the ``words are tools"\footnote{``Think of the tools in a tool-box: there is a hammer, pliers, a saw, a screw-driver, a rule, a glue-pot, glue, nails and screws. -- The functions of words are as diverse as the functions of these objects. (And in both cases there are similarities.)", \S 11, Part I, p. 6$^\mathrm{e}$ of the \emph{Investigations}.

Consider this passage in an unpublished manuscript of 1929--30 (Ms.\ 107):
\begin{quote}
Kann man sagen: der Sinn des || eines Satzes ist sein Zweck? [Oder von einem Wort ,,its meaning is its purpose".]
\end{quote}
which may be translated as:
\begin{quote}
Can one say: the sense of the || of a sentence is its purpose? [Or  a word ``its meaning is its purpose".]
\end{quote}
Ms.\ 107,234 (Wittgenstein Nachlass Ms-107: III, Philosophische Betrachtungen (\"ONB, Cod.Ser.n.22.020) - (IDP). (Accessed 11 Jun 2022))} conceptual link made mainly in Part I of the \emph{Investigations}, but turns up in the beginning of the later period (Ms.\ 109, 1930--31) in the form of `meaning as explanation' and `meaning as expectation',\footnote{Wittgenstein states:
\begin{quote}
Die Realit\"at ist keine Eigenschaft die dem Erwarteten noch fehlt \& die nun hinzutritt wenn es eintritt. -- Sie ist auch nicht wie das Tageslicht das den Dingen erst ihre Farbe gibt wenn sie -- im Dunkeln schon gleichsam farblos vorhanden sind.
     Alle diese grammatischen Formen stellen den Gegensatz Erwartung \& Erf{\"u}llung nicht dar. Die Sprache stellt ihn nur so dar, wie sie ihn immer darstellt durch den Gegensatz der S\"atze ,,ich erwarte p" \& ,,p ist eingetroffen".
	 			 	
Wie konnte ich es erwarten, \& es kommt dann wirklich; -- als ob die Erwartung ein dunkles Transparent w{\"a}re \& mit der Erf{\"u}llung
das Licht dahinter angez{\"u}ndet w{\"u}rde. Aber jedes solche Gleichnis ist falsch weil es die Realit{\"a}t als einen beschreibbaren Zusatz zur Erwartung || zum Gedanken darstellt, was unsinnig ist.
\end{quote}
which may be translated as:
\begin{quote}
Reality isn't a quality that's missing from what's expected \& that's added when it happens. -- Nor is it like daylight, which only gives things their color when in the dark they are, as it were, colorless.
     None of these grammatical forms represent the antithesis of expectation \& fulfillment. Language only represents it as it always represents it through the opposition of the sentences ``I expect p" \& ``p has arrived".

How could I expect it, \& then it really comes; -- as if the expectation were a dark transparency \& with the fulfillment
the light behind would be lit. However,
any such simile  is false because it presents reality as a describable addition to
expectation, to thought, which is nonsense.
\end{quote}
Ms-109,54--55, (Wittgenstein Nachlass Ms-109: V, Bemerkungen (WL) - (IDP). (Accessed 11 Jun 2022))
}
 seems to find its roots in very early remarks documented in early writings such as the \emph{Notebooks} and the \emph{Tractatus}. If the former shows the early signs of the connections between meaning and use -- ``The way in which language signifies is mirrored in its use"\footnote{dated 11.9.16, p.\ 82$^\mathrm{e}$ of the \emph{Notebooks 1914-1916}.}-- in the latter (and in its early version, the so-called \emph{Prototractatus}) one can find remarks such as:
\begin{quote}
In order to recognize the symbol in the sign one pays attention to the use.\footnote{\S 3.252, \emph{Prototractatus}. Wittgenstein \emph{Nachlass} Ms-104: Logisch-Philosophische Abhandlung [so-called \emph{Prototractatus}] (Bodl, MS. German d. 7) - (IDP) (accessed 13 May 2021)

In fact, this remark is already present in an earlier manuscript, i.e.\ Ms 101, dated 23.10.14: 
\begin{quote}
Um das Zeichen im Zeichen zu erkennen mu\ss\ man auf den Gebrauch achten. (Ms-101,59r)
\end{quote}
which may be translated as: 
\begin{quote}
In order to recognise the sign in the sign, one must pay attention to the use.
\end{quote}}
\end{quote}
\begin{quote}
In der Philosophie f\"uhrt die Frage ``wozu gebrauchen wir eigentlich jenes Wort, jenen Satz" immer wieder zu wertvollen Einsichten. [In philosophy, the question ``to what end do we actually use that word, that sentence" repeatedly leads to valuable insights.]\footnote{\S 6.1221, \emph{Prototractatus}. Wittgenstein \emph{Nachlass} Ms-104: Logisch-Philosophische Abhandlung [so-called \emph{Prototractatus}] (Bodl, MS. German d. 7) - (IDP) . (accessed 15 May 2022)}
\end{quote}
\begin{quote}
If a sign is useless, it is meaningless. That is the point of Occam's maxim.\footnote{\S 3.328, p.\ 16 of the \emph{Tractatus}.  This is already present in a 1915 manuscript:
\begin{quote}

5'30634\qquad Zeichen die Einen Zweck erf\"ullen sind logisch \"aquivalent. Zeichen die keinen Zweck erf\"ullen logisch bedeutungslos.
 	 			 	
3'2521\qquad Wird ein Zeichen nicht gebraucht, so ist es bedeutungslos. Das ist der Sinn der Devise Ockhams.
\end{quote}
which may be translated as:
\begin{quote}
5'30634\qquad Signs that serve one purpose are logically equivalent. Signs that serve no purpose are logically meaningless.
 
3'2521\qquad If a sign is not used, it is meaningless. That is the point of Occam's motto.
\end{quote}
Ms-104,60 (1915--18) (In: Wittgenstein, Ludwig: Interactive Dynamic Presentation (IDP).) (Accessed 16 Jun 2022)}
\end{quote}
and `useless' here is to be understood as `purposeless', as it becomes clear later in the same \emph{Tractatus}:
\begin{quote}
Occam's maxim is, of course, not an arbitrary rule, not one that is justified by its success in practice: its point is that unnecessary units in a sign-language mean nothing.
Signs that serve one purpose are logically equivalent, and signs that serve none are logically meaningless.\footnote{Ibid., \S 5.47321, p.\ 48.}
\end{quote}

A remark dated 23.4.15, once again makes a connection between meaning and use/purpose/ usefulness:
\begin{quote}
Es ist klar da\ss\ Zeichen, die denselben Zweck erf\"ullen logisch identisch sind. Das rein Logische ist eben das was alle diese leisten k\"onnen.\footnote{which may be translated as:
\begin{quote}
It is clear that signs serving the same purpose are logically identical. The purely logical is what all these can achieve.
\end{quote}
Ms-102,75r (Ms-102 (WL)  (IDP).) (Accessed 16 Jun 2022.)
}
\end{quote}

For the former, the question of whether the logical constants exist is not even considered relevant, since they can ``disappear". And again the position can be seen as coherent with that expressed in the remark above, considering such ``constants" as implements, tools that make sense in their use/usefulness. In this particular respect, Wittgenstein's words are phrased in the {\em Notebooks} as follows:
\begin{quote}
With the logical constants one need never ask whether they exist, for they can even vanish!\footnote{dated 25.10.14, p.\ 19$^\mathrm{e}$ of the \emph{Notebooks}.}
\end{quote}

Further down in the \emph{Notebooks} the (introduction of the) sign `0' is said to get its meaning from its usefulness in making the decimal notation possible:

\begin{quote}
The introduction of the sign ``0" in order to make the decimal notation possible: the logical significance of this procedure.\footnote{dated 16.11.14, p.\ 31$^\mathrm{e}$ of the \emph{Notebooks}}
\end{quote}

If meaning is to be understood as use/usefulness, logical connectives are effectively instruments, rather than constants.This Wittgenstein declares, not only in later writings such as:
\begin{quote}
Es hat keinen Sinn von S\"atzen zu reden die als Instrumente keinen Wert haben.

Der Sinn eines Satzes ist sein Zweck.\footnote{which may be translated as:
\begin{quote}
There is no point in talking about sentences that have no value as instruments.

The sense of a sentence is its purpose.
\end{quote}
Ms-107,249 (1929--30) (Wittgenstein \emph{Nachlass} Ms-107: III, Philosophische Betrachtungen (\"ONB, Cod.Ser.n.22.020) - (IDP).)
}
\end{quote}
even in the \emph{Tractatus} where the truth-values based semantical view inherited from Frege and Russell characterises the `early Wittgenstein' as distinct from the `later Wittgenstein and the use-based semantical view'. Such an original perspective on logic and language versus reality already appears when he says that logical connectives are tools, therefore the idea that `logical connectives are constants' and the viewpoint that `there are logical constants, namely ``$\lor$", ``$\supset$", ``$.$",  (``$\land$"), etc.' are hardly tenable. In order to make the point he recalls the usual interdefinability of such `constants':\footnote{When talking of logical `constants' one is concerned to distinguish these from the variables on which they operate in propositions. In this context, the crucial point is Frege's or Russell's distinction between  `logical objects' and `logical constants'. (Thanks to an anonymous reviewer for such an important reminder.)}
\begin{quote}
At this point it becomes manifest that there are no `logical objects' or `logical constants' (in Frege's and Russell's sense).\footnote{\S 5.4, p.\ 44 of the \emph{Tractatus}.}
\end{quote}

At this point it would not be unreasonable to say that insofar as being understood as implements, instruments, the logical connectives would not have to be looked at as `primitive' signs. And  Wittgenstein insists on this point, already in the \emph{Tractatus}:
\begin{quote}
It is self-evident that $\lor$, $\supset$, etc. are not relations in the sense in which right and left etc.\ are relations.
The interdefinability of Frege's and Russell's `primitive signs' of logic is enough to show that they are not primitive signs, still less signs for relations.\footnote{Ibid., \S 5.42, p.\ 44.}
\end{quote}

Thus, despite resorting to truth-values and truth-functions when building a precise account of meaning for logical propositions by
``a proposition is a truth-function of elementary propositions. (an elementary proposition is a truth-function of itself.),"\footnote{Ibid., \S 5, p.\ 36.} Wittgenstein did not sound convinced that the concept of truth-function could work equally well for all logical concepts and connectives of Frege-Russell's logical calculus. As a sort of `canonical' example, to which he returns in later writings, the generality concept in which the universal quantifier is embedded was explicitly dissociated from truth-functions in:
\begin{quote}
I dissociate the concept all from truth-functions.

Frege and Russell introduced generality in association with logical product or logical sum. This made it difficult to understand the propositions `$(\forall x).fx$' and `$(\exists x).fx$', in which both ideas are embedded.\footnote{Ibid., \S 5.521, p.\ 51.}
\end{quote}

 The use of `dissociate' in the translation by D.\ Pears and B.\ McGuinness appears to reinforce Wittgenstein's deliberate intention in ``Ich trenne den Begriff Alle von der Wahrheitsfunktion" to disconnect one from another. This is clear from the correspondence between Wittgenstein and the first (English) translator C.\ K.\ Ogden. In the `Comments II on Separate Sheets', from \emph{Letters to C. K. Ogden}, p.\ 61, he says:
\begin{quote}
5.521 I mean ``separate" and not ``derive". They were connected and I separate them.
\end{quote}
and the comments by the Editor, p.\ 65:
\begin{quote}
5.521 The passage under discussion runs: ``I separate the concept all from the truth-function."

Ogden had queried whether ``separate" should not be changed to ``derive".
\end{quote}

Later in a typescript of 1931--32 (partly published in \emph{Philosophical Grammar} and \emph{Remarks on Frazer's Golden Bough}) one can find again a reference to the fact that the existential quantifier should not be viewed as a mere sum:
\begin{quote}
Ja sagt denn eben $(\exists x).fx\lor fa = (\exists x).fx$ nicht, da\ss\ $fa$ schon in $(\exists x).fx$ enthalten ist? Zeigt es nicht die Abh{\"a}ngigkeit des $fa$ vom $(\exists x).fx$? Nein, au{\ss}er, wenn $(\exists x).fx$ als logische Summe definiert ist (mit einem Summanden $fa$).  Ist das der Fall, so ist $(\exists x).fx$ (nichts als) eine Abk{\"u}rzung.\footnote{which may be translated as:
\begin{quote}
Doesn't $(\exists x).fx\lor fa = (\exists x).fx$ mean that $fa$ is already contained in $(\exists x).fx$ ? Doesn't it show the dependence of $fa$ on $(\exists x).fx$ ? No, except when $(\exists x).fx$  is defined as a logical sum (with a summand $fa$).  If this is the case, then $(\exists x).fx$  is (nothing but) an abbreviation.
\end{quote}
Ts 211,38 (Wittgenstein \emph{Nachlass} Ts-211 (WL) (IDP).) (Accessed 18 June 2022)
}

\end{quote}

From this it should be clear that it is unwarranted to say that Wittgenstein in the \emph{Tractatus} ``construes quantified sentences as conjunctions and disjunctions"\footnote{p.\ 110 of (Hintikka \& Hintikka (1986)).}. Hintikka \& Hintikka's conclusion seems to draw on what Moore reports.\footnote{Moore says that Wittgenstein in 1930--33, ``went on to say that one great mistake which he made in the \emph{Tractatus} was that of supposing that in the case of all classes ``defined by grammar", general propositions were identical with logical products or with logical sums (meaning by this logical products or sums of the propositions which are values of $fx$) as, according to him, they really are in the case of the class of ``primary colours". He said that, when he wrote the \emph{Tractatus}, he had supposed that all such general propositions were ``truth-functions"; (...) He said that, when he wrote the \emph{Tractatus}, he would have defended the mistaken view which he then took by asking the question: How can $(x).f x$ possibly entail $f a$, if $(x).f x$ is not a logical product? And he said that the answer to this question is that where $(x).fx$ is not a logical product, the proposition ``$(x).fx$ entails $fa$" is ``taken as a primary proposition", whereas it is a logical product this proposition is deduced from other primary propositions.", Moore (1955, 3).} Wittgenstein might well have said that, given that he was very critical, almost ruthless, with himself and he often insisted on the big mistakes he made in the \emph{Tractatus}. For example, he says in a section entitled `Criticism of my former view of generality' of the \emph{Philosophical Grammar}, p.\ 268, that his ``view about general propositions was that $(\exists x).\varphi(x)$  is a logical sum". Nonetheless, nowhere in the \emph{Tractatus} is the existential quantifier treated as a logical sum. Moreover, as we shall see below, early remarks bear witness to the different, `non-truth-functional' treatment he gave to the universal quantifier.

For the particular case of the logical notion of generality and the corresponding symbol for the universal quantifier, since the early days of Wittgenstein's own writings, the notion of `following from' has prevailed to the more \emph{Tractatus}-characteristic truth-function. At some point, notably in the bulk of the \emph{Tractatus}, Wittgenstein seemed to have preferred the more \emph{Tractatus}-characteristic truth-function, but even in this book one can see traces of the more prevalent across essentially both early as well as late writings.:
\begin{quote}
(24.11.14.) Da\ss\ man aus dem Satz ,,$(x) \phi x$" auf den Satz ,,$\phi a$" schlie{\ss}en kann das zeigt wie die Allgemeinheit auch im Zeichen ,,$(x) \phi x$" vorhanden ist.

Und das gleiche gilt nat\"urlich f\"ur die Allgemeinheitsbezeichnung \"uberhaupt.\footnote{which may be translated as:
\begin{quote}
(24.11.14.) The fact that one can deduce the sentence ``$\phi a$" from the sentence ``$(x) \phi x$" shows how generality is also present in the sign ``$(x) \phi x$".

And the same applies, of course, to the generality designation as such.
\end{quote}
Ms-102,38r (Wittgenstein \emph{Nachlass} Ms-102 (WL) (IDP). (Accessed 11 Jun 2022) )
}
\end{quote}
A similar remark also appeared in the \emph{Tractatus}:
\begin{quote}
 The possibility of inference from $(x).fx$ to $fa$ shows that the symbol $(x).fx$ itself has generality in it.\footnote{\emph{Tractatus}, \S 5.1311, p.\ 38.}
 \end{quote}
 
The inference demonstrating the (immediate) consequences one can draw from the proposition `$(x).fx$' actually makes clear that the distinguishing characteristic of the major logical sign (namely `$(x).$', \emph{for all} $x$) is that of generality. From this standpoint, one can read many of the remarks contained in the \emph{Tractatus} itself in a way which has, to the best of our knowledge, hardly been done in the literature on the early Wittgenstein, and place them still in line with the remarks from the later Wittgenstein using the same connections between meaning, use/usefulness and consequences. When explaining the consequences one can draw from an existential statement, Wittgenstein must resort to other resources, obviously, which is reflected in the indirect form of drawing consequences from $(\exists x)f.(x)$ used in the rule of $\exists$-elimination in Natural Deduction.\footnote{In TS 212 there is a long discussion of the existential quantifier and the explanation of its meaning; here is a relevant excerpt:
\begin{quote}
$(\exists x).fx\lor fa = (\exists x).fx, \ (\exists x).fx\ \&\ fa = fa$ Wie wei\ss\ ich das? (denn das Obere habe ich sozusagen bewiesen). Man m\"ochte etwa sagen: ``ich verstehe `$(\exists x).fx$' eben". (Ein herrliches Beispiel dessen, was `verstehen' hei{\ss}t.)

     Ich k\"onnte aber ebensogut fragen ``wie wei\ss\ ich, da\ss\ $(\exists x).fx$ aus $fa$ folgt" und antworten: ``weil ich `$(\exists x).fx?$' verstehe". Wie wei\ss\ ich aber wirklich, da\ss\ es folgt? -- Weil ich so kalkuliere.
 
 Wie wei\ss\ ich, da\ss\ $(\exists x).fx$ aus $fa$ folgt? Sehe ich quasi hinter das Zeichen ``$(\exists x).fx$", und sehe den Sinn, der hinter ihm steht und daraus || aus ihm, da\ss\ er aus $fa$ folgt? ist das das Verstehen?
 
     Nein, jene Gleichung ist ein Teil des Verstehens || Verst\"andnisses || dr\"uckt einen Teil des Verstehens || Verst\"andnisses aus (das so ausgebreitet vor mir liegt).
     (Denke an die  || Vergleiche die Auffassung des Verstehens, das urspr\"unglich mit einem Schlag erfa{\ss}bar || ein Erfassen mit einem Schlag, erst so ausgebreitet werden kann.
     Wenn ich sage ``ich wei{\ss}, da\ss\ $(\exists x).fx$ folgt, weil ich es verstehe", so hie{\ss}e das, da\ss\ ich, es verstehend, etwas Anderes sehe, als das gegebene Zeichen, gleichsam eine Definition des Zeichens, aus der das Folgen hervorgeht.
\end{quote}
which may be translated as:
\begin{quote}
$( \exists x).fx \lor fa = ( \exists x).fx$, $( \exists x).fx \& fa = fa$ How do I know? (for I have proved the above, so to speak). One would like to say something like: ``I understand `$( \exists x).fx$' just". (A splendid example of what `understand' is.)

But I might as well ask ``how do I know that $( \exists x).fx$ follows from $fa$" and reply: ``because I understand `$( \exists x).fx$'?". But how do I really know that it follows? - Because I calculate so.

How do I know that $( \exists x).f x$ follows from $fa$? Do I see behind the sign ``$( \exists x).f x$"'', and see the sense that stands behind it and out of it  that it follows from $fa$? is that understanding?

No, that equation is a part of understanding || understandings || a part of understanding || understandings (which thus displayed before me). (Think of them? || Compare the notion of understanding that can be grasped in one stroke || capturing in one fell swoop that can only be spread out. When I say, ``I know that `$( \exists x).fx$' follows because I understand it," this means that, understanding it, I see something other than the given character, as it were, a definition of the character from which the consequences emerge.
\end{quote}
Ts-212,IX-66-2 (Wittgenstein \emph{Nachlass} Ts-212 (WL) - (IDP). (with Heinz Wilhelm Kr\"{u}ger (2022: proofreading of normalized transcription))) (Accessed 12 June 2022).
}
For example, a remark concerning the view that in logic the performance of deductions (`\emph{Folgern}') occurs a priori, and it appears that one can understand this in the light of, e.g., Hintikka's semantical games, a (semi-)formal counterpart to later Wittgenstein's `the meaning of a sign is its use in the language-game', where the logical connectives are given meaning by rules of `one-step' deduction. And indeed, if ``before a proposition can have a sense, it must be completely settled what propositions follow from it"\footnote{\S 3.20103, p.\ 65 of the \emph{Prototractatus}. An almost identical remark already appears in the \emph{Notebooks 1914--16}: ``We might demand definiteness in this way too: if a proposition is to make sense then the syntactical employment of each of its parts must be settled in advance. -- It is, e.g., not possible only subsequently to come upon the fact that a proposition follows from it. But, e.g., what propositions follow from a proposition must be completely settled before that proposition can have a sense!", \emph{Notebooks}, dated 18.6.15, p.\ 64$^\mathrm{e}$.}, according to a remark dated 15.10.14 in Ms101:
\begin{quote}
Man kann sagen: die v\"ollig allgemeinen S\"atze kann man alle a priori bilden.\footnote{which may be translated as:
\begin{quote}
One can say: the completely general propositions can all be formed a priori.
\end{quote}
Ms-101,43r (Wittgenstein \emph{Nachlass} Ms-101 (WL)
(IDP).) (Accessed 19 June 2022).}
\end{quote}
as well as the remark ``5'0415 Alles Folgern geschieht a priori"\footnote{Wittgenstein \emph{Nachlass} Ms-104: Logisch-Philosophische Abhandlung [so-called \emph{Prototractatus}] (Bodl, MS. German d. 7)
(IDP).
}
 of Ms-104 (published as \emph{Prototractatus}), so that:
\begin{quote}
All deductions are made a priori.\footnote{\S 5.133, p.\ 39, of the \emph{Tractatus}.}
\end{quote}
according to the \emph{Tractatus}. It sounds appropriate to say these remarks can be taken along the same lines advocated here by simply acknowledging that to show the (immediate) consequences which can be drawn from a proposition is to demonstrate what can be done with it, those a priori, immediate deductions belonging to the `definition'. Moreover, in the same line of reasoning one has to say that (using early Wittgenstein's own words), ``logic must look after itself"\footnote{Ibid., \S 5.473, p.\ 47.} simply because, for a proposition to be defined and have a sense in logic, the rule showing the consequences (deductions) one can draw from the proposition must be given a priori. Therefore, ``we cannot give a sign the wrong sense"\footnote{Id. Ibid., \S 5.4732.}, and so, the (one-step) `deductions' are made a priori just because they constitute the meaning of the proposition.

At this point, it makes sense to once again emphasise the remarks suggesting the links between the meaning of an assertion and what can be deduced from it, even in very early writings, such as: ``A proposition affirms every proposition that follows from it."\footnote{Ibid., \S 5.124, p.\ 38.}, and: `` `$p.q$' is one of the propositions that affirm `$p$' and at the same time one of the propositions that affirm `$q$'.\footnote{Id. Ibid., \S 5.1241.} of the \emph{Tractatus}.
As early as 1912, in a letter to Russell, dated 1.7.12, he had already expressed
something which is essentially similar to what is referred to above:
\begin{quote}
Will you think that I have gone mad if I make the following suggestion?: The sign $(x).\phi x$ is not a complete symbol but has meaning only in an inference of the kind: from $\vdash\phi x\supset_x\psi x.\phi(a)$ follows $\psi(a)$. Or more generally: from $\vdash (x).\phi x.\varepsilon_0(a)$ follows $\phi(a)$. I am -- of course -- most uncertain about the matter but something of the sort might really be true.\footnote{Letter R.3, p.\ 12 of Wittgenstein (1974) \emph{Letters to Russell, Keynes and Moore}}
\end{quote}
The universal quantifier over $\phi(x)$ has meaning when one explains that the immediate consequence of the proposition is an instance of the sentence for a generic element $a$, i.e., $\phi(a)$.

Studies based on the catalogued \emph{Nachlass} have uncovered many misconceptions generated by the rather fragmented character of Wittgenstein's philosophical publications. In the light of those findings, it seems reasonable to argue, for example, that the logical atomism of the \emph{Tractatus} was rejected by Wittgenstein in so far as it assumed that a name denoted an object, and the object was its `meaning' (i.e.\ denotation). The opening lines of the \emph{Investigations} (and beginning of: \emph{MS 111, TS 213, Blue Book, Brown Book}; cf.\ also Ms.\ 109 (1930--31) \emph{Bemerkungen}\footnote{Ms.109,140:
\begin{quote}
Die Bedeutung (eines Wortes) k\"{o}nnte || kann nur das sein was wir in der Erkl\"{a}rung eines Wortes erkl\"{a}ren.
 			 	
Man hat immer die falsche Vorstellung als w\"{a}re die Bedeutung eines Wortes ein || handelte es sich bei der Bedeutung eines Wortes um einen Gegenstand d.h.\ ein Ding, in dem Sinn in dem das Schwert Nothung die Bedeutung des Wortes ,,Nothung" war. Aber auch hier stimmt etwas nicht, denn ich kann doch sagen ,,das Schwert Nothung || ,,Nothung existiert nicht mehr" \& ist etwa hier ,,Nothung" bedeutungslos eben weil das Schwert nicht mehr existiert?
\end{quote}
which may be translated as:
\begin{quote}
The meaning (of a word) could be || can only be what we explain in the explanation of a word.
 
One always has the misconception that the meaning of a word is || the meaning of a word were an object i.e. a thing, in the sense in which the sword Nothung was the meaning of the word ``Nothung". But there is something wrong here too, because I can still say ``the sword Nothung || 'Nothung' no longer exists" \& is 'Nothung' meaningless here precisely because the sword no longer exists?
\end{quote}
(Wittgenstein \emph{Nachlass} Ms-109: V, Bemerkungen (WL) - (IDP).)}) reveal one of the main criticisms against the \emph{Tractatus}' atomism: the assumption that names denote objects which are their meanings. In the latter phase, whilst that assumption is dropped,  atomism is retained, leaving a sort of `functional atomism', where words are taken to be individually identifiable devices with a meaning (function, purpose, usefulness) in the calculus of language. A revision of his early (`official') strict atomism into a kind of manageable holism accounts for both the continuity (i.e.\ atomism, with a `functional' (`manageably holistic') ingredient, as David Pears emphasises: ``the idea that a sentence can only be explained by their places in systems (...) is one of the central points of the picture theory"\footnote{footnote 34, p.\ 121 of Pears (1987).} and ``Everyone is aware of the holistic character of Wittgenstein's later philosophy, but it is not so well known that it was already beginning to establish itself in the \emph{Tractatus}."\footnote{Opening lines of Pears (1990).}) and the discontinuity (i.e.\ `meaning as denotation' vs.\  `meaning as use') of his core views.

The challenge now is to build an adequate basis on which a functional atomism can be used to construct a framework where the logical connectives have a meaning (function, purpose, usefulness) in the calculus of language which is made clear by the explanation of the (immediate) consequences one is supposed to draw from the corresponding proposition.We have discussed elsewhere (2001, 2008) how the so-called `functional interpretation' of logics can help. As an additional step towards building this framework, here we show connections with the ``pragmatist"/``dialogical" approaches to meaning, which was already in Peirce's works on the   \emph{Interpreter} vs \emph{Utterer} dichotomy, improve on the view of logical constants as tools whose meanings are determined by their use/usefulness. Augmenting those previous publications, this connection appears to hinge upon the ``game"/``dialogue" approaches to meaning (Lorenzen, Hintikka, Fra\"{\i}ss\'e). 
The (Curry--Howard) functional interpretation of logical connectives is an `enriched' system of Natural Deduction, terms representing proof-constructions are carried alongside the formulas. It can be used to formulate various logical calculi operationally, where the way one handles assumptions determines the logic. So, the novelty of pointing to the connection with Peirce's dichotomy reinforces the fact that the proof reduction rules in the Curry-Howard  approach reveals that  it is naturally atomistic, still concerned with proofs, and is functional in the sense that the main semantical device is that of reduction of proofs (not only assertability conditions), finding a parallel in the dialogue/game-like semantics of Lorenzen (and Hintikka): canonical terms (introduction) stand for the `Assertion' (resp.\ `Myself') moves, and non-canonical terms (elimination) stand for the `Attack' (resp.\ `Nature') moves. This appears to accord with Wittgenstein's `functional atomism': it accounts for the early atomism (already with holistic tendencies according to D.\ Pears (1990)) in a functional fashion, abiding by `words are tools' and `language is a calculus' (which, S.\ Hilmy verifies, ``was already in place in early 1930s, and remained at least until 1946").\footnote{Hilmy (1987, p.\ 14) points out that even in the \emph{Investigations} one finds the use of the calculus/game paradigm to  understand language, such as ``im Laufe des Kalk\"uls" (Part I, \S 559) (in the English translation: ``in operating with the word"), and ``it plays a different part in the calculus". Hilmy also quotes from a late (1946) unpublished manuscript (MS 130), ``this sentence has use in the calculus of language" is compatible with ``asking whether and how a proposition can be verified is only a particular way of asking `How do you mean?' " (\emph{Investigations}, Part I, \S 353).}

In the writings of the early period of his later phase, documented in the notebooks containing remarks mainly on the foundations of mathematics now published as \emph{Philosophical Remarks} and \emph{Philosophical Grammar}, there is a clear shift in focus from truth-values and truth-functions, characteristic of the \emph{Tractatus}, to mathematical proofs. The concept of mathematical proof is thoroughly investigated, and in some passages it is even placed in a rather privileged position with respect to the semantics of mathematical language. One of the lessons taught by those `transitional' writings seems to be that (using Wittgenstein's own words) ``when one wants to know the meaning, the sense, of a proposition, one must look at what its proof proves". Here is a remark from Ms 106: II (p.\ 182) (1929):
\begin{quote}
Man k{\"o}nnte fragen: Was sagt (x)2x = x + x? Es sagt da\ss\ alle Gleichungen von der Form 2x = x + x richtig sind. Aber hei{\ss}t das etwas? Kann man sagen: Ja ich sehe da\ss\ alle Gleichungen dieser Form richtig sind, so kann ich jetzt schreiben ,,(x)2x = x + x"?
	 			 	
Ihre Bedeutung mu\ss\ aus ihrem Beweis hervorgehen. Was der Beweis beweist das ist die Bedeutung des Satzes (nicht mehr \& nicht weniger).\footnote{which may be translated as:
\begin{quote}
One might ask: What does (x)2x = x + x say? It says that all equations of the form 2x = x + x are correct. But does that mean something? Can one say: Yes, I see that all equations of this form are correct, so I can now write ``(x)2x = x + x"?
 
Their meaning must be evident from their proof. What the proof proves is the meaning of the proposition (no more \& no less).
\end{quote}
Wittgenstein \emph{Nachlass} Ms-106: II (\"ONB, Cod.Ser.n.22.019) - (IDP). (accessed 16 May 2022).}
\end{quote}

It does not take too long, however, for an important revision of the position on the connections between meaning and proofs to be made; in the \emph{Remarks on the Foundations of Mathematics}, written between 1941--44:
\begin{quote}
I once said: `If you want to know what a mathematical proposition says, look at what its proof proves'.\footnote{Cf.\ \emph{Philosophical Grammar}, p.\ 369; \emph{Philosophical Remarks}, pp.\ 183--4.} Now is there not both truth and falsehood in this? For is the sense, the point, of a mathematical proposition really clear as soon as we can follow the proof?\footnote{\S 10, Part VII, p.\ 367 of the \emph{Remarks on the Foundations of Mathematics}.}
 \end{quote}
A fair number of remarks which reveal a much less explicit statement of this revision regarding the connections between the meaning of a proposition and its proof appear in later stages of the later period. One of them, from the same \emph{Remarks on the Foundations of Mathematics} written around 1943/1944, states:
\begin{quote}
The proof of a proposition certainly does not mention, certainly does not describe, the whole system of calculation that stands behind the proposition and gives it its sense.\footnote{Ibid., \S 11, Part VI, p.\ 313.}
\end{quote}
It reveals his concern about the insufficiency of explaining what counts as a proof (i.e. explanation of assertability conditions) with respect to explaining meaning.\footnote{Consider here this remarks written in 1929 in Ms-106:
\begin{quote}
The mathematical proposition then relates to its proof as the uppermost surface of a body relates to itself. One could speak of the proof body of the proposition.
 
Only if there is a body behind the surface has the sentence meaning for us.
 
One also says: the mathematical proposition is (only) the last link in a chain of proof. (Ms 106, p.\ 195)
\end{quote}
Wittgenstein \emph{Nachlass} Ms-106: II (\"ONB, Cod.Ser.n.22.019) (IDP) (accessed 13 May 2021)
}
 The apparently unworkable holism is shown perfectly feasible to the formulation of a theory of meaning for the language of mathematics by his suggestions on the connections between the meaning of a proposition and the explanation of the immediate consequences one can draw from it:
\begin{quote}
One learns the meaning [Bedeutung] of `all' by learning that `$fa$' follows from `$(x).fx$'.

(...)

And the meaning [Sinn] of `$(x).fx$' is made clear by our insisting that
`$fa$' follows from it.\footnote{Ibid., Part I, \S 10, p.\ 41 and \S 11, p.\ 42, respectively.}
\end{quote}
These parts of language are understood only in relation to other
neighbouring parts, as D.\ Pears (1987) points out: 
\begin{quote}The new idea [of W.'s later theory of language] was that we experience each part of language in its relations to other, neighbouring parts. (p.\ 171)
\end{quote}
Wittgenstein's use of the idea that ``the meaning of a word is determined by its r\^ole in the language" would be understood in the following way: by showing how to `eliminate', `get rid of' the logical sign, demonstrating the conclusions one can draw from a proposition which has it as its major connective, one  demonstrates the r\^ole it plays in the calculus/system of propositions.

In Ms-145 [so-called C1]  (1933) one finds remarks to the effect of associating the meaning of a sentence with the consequences one can draw from it:
\begin{quote}
Der Satz ist f\"ur mich kein blo{\ss}er Laut, er ruft in mir eine Vorstellung hervor etwa von jemandem in dem Senathaus.
Aber der Satz \& diese Vorstellung ist nicht blo\ss\ ein Laut \& eine schwache Vorstellung sondern der Satz hat es (sozusagen) in sich diese Vorstellung hervorzurufen aber auch andere Konsequenzen \& das ist sein Sinn. Die Vorstellung scheint nur ein schwaches Abbild dieses Sinnes oder, sagen wir, nur eine Ansicht dieses Sinnes. -- Aber was meine ich denn damit; sehe ich den Satz eben nicht als Glied in einem System von Konsequenzen? Kann ich nicht sagen: Wenn ich mir {\"u}berlege, was war der Sinn eines Satzes  || worin es bestehe da\ss\ der Satz einen Sinn habe, so finde ich da{\ss} es || er darin liegt da\ss\ so aus dem Satz nach Festsetzungen Konsequenzen hervorgehen, da\ss\ ich den Satz {\"a}hnlich wie den Zug in einem Spiel betrachte.\footnote{which may be translated as:
\begin{quote}
The phrase isn't just a sound to me, it conjures up an image in my mind, as of someone in the Senate House.
But the sentence \& this idea is not just a sound \& a weak idea. The sentence has it (so to speak) in itself to evoke this idea but also other consequences \& that is its sense. The idea seems only a faint image of that sense, or, shall we say, only a view of that sense. -- But what do I mean by that; don't I see the sentence as a link in a system of consequences? Can't I say: If I think about it, what was the sense of a sentence  || in what it consists the proposition a sense, I find that it || it lies in the fact that consequences emerge from the proposition according to stipulations, that I regard the proposition in a way similar to the move in a game.
\end{quote}
Ms-145,80-81 (1933)  (In: Wittgenstein, Ludwig: Interactive Dynamic Presentation (IDP).) (Accessed 13 June 2022)
}
\end{quote}

Wittgenstein's acknowledging the inappropriateness of using mathematical proofs and assertability conditions as semantical notions appears to have been neglected by some `verificationist' theories of meaning which seek a conceptual basis in Wittgenstein's semantic theory and suggest that `meaning is use' is their underlying semantical paradigm. For example, Michael Dummett (1977, 12):
\begin{quote}
the meaning of each [logical] constant is to be given by specifying, for any sentence in which that constant is the main operator, what is to count as a proof of that sentence, it being assumed that we already know what is to count as a proof of any of the constituents.\footnote{As a side remark, what we have here may be understood as showing that constructivist semantics need not and should not be `verificationist', in (at least) Dummett's sense. The point to be retained from (e.g.) Brouwer's intuitionism is its constructivism. The 'solipsism' is spurious; but, in this sense, Brouwer and Wittgenstein are correct that understanding must be 1st person, though anyone's achieving such understanding can of course be within social and worldly contexts, including chalk boards. (Thanks to anonymous reviewer for pointing this out.)}
\end{quote}
Underlying the so-called Dummett--Prawitz theory of meaning via proofs, D.\ Prawitz (1977), puts forward an interpretation of Wittgenstein's `meaning is use':
\begin{quote}
As pointed out by Dummett, this whole way of arguing with its stress on communication and the role of the language of mathematics is inspired by ideas of Wittgenstein and is very different from Brouwer's rather solipsistic view of mathematics as a languageless activity. Nevertheless, as it seems, it constitutes the best possible argument for some of Brouwer's conclusions. (...)

I have furthermore argued that the rejection of the Platonistic theory of meaning depends, in order to be conclusive, on the development of an adequate theory of meaning along the lines suggested in the above discussion of the principles concerning meaning and use. Even if such a Wittgensteinian theory did not lead to the rejection of classical logic, it would be of great interest in itself.
\end{quote}
Along the same lines, Dummett (1991, 251f) returns to it:
\begin{quote}
Gerhard Gentzen, who, by inventing both natural deduction and the sequent calculus, first taught us how logic should be formalised, gave a hint how to do this, remarking without elaboration that `an introduction rule gives, so to say, a definition of the constant in question', by which he meant that it fixes its meaning, and
that the elimination rule is, in the final analysis, no more than a consequence of this definition. (...) Plainly,
the elimination rules are not consequences of the introduction rules in the straightforward sense of being derivable from them; Gentzen must therefore have
had in mind some more powerful means of drawing consequences. He was also implicitly claiming that
the introduction rules are, in our terminology, self-justifying.
\end{quote}
In the sequel, Dummett suggests that P.\ Martin-L\"of in his intuitionistic type theory would be advocating a different standpoint, namely, that the rules of elimination, which establish the immediate consequences of a proposition, define the meaning of the logical connective:
\begin{quote}
Intuitively, Gentzen's suggestion that the introduction rules be viewed as fixing the meanings of the logical constants has no more force than the converse suggestion, that they are fixed by the elimination rules; intuitive plausibility oscillates between these opposing suggestions as we move from one logical constant to another. Per Martin-L\"of has, indeed, constructed an entire meaning-theory for the language of mathematics on the basis of the assumption that it is the elimination rules that determine meaning. The underlying idea is that the content of a statement is what you can do with it if you accept it -- what difference learning that it is true will, or at least may, make to you. This is, of course, the guiding idea of a pragmatist meaning-theory. When applied to the logical constants, the immediate consequences of any logically complex statement are taken as drawn by means of an application of one of the relevant elimination rules.
\end{quote}
Assuming a similar reasoning, Martin-L\"of himself clearly adheres to the same view that the introduction rules, i.e.\ those which establish the assertability conditions, should be seen as the rules defining the meaning of a logical connective:
\begin{quote}
The intuitionists explain the notion of proposition, not by saying that a proposition is the expression of its truth conditions, but rather by saying, in Heyting's words, that a proposition expresses an expectation or an intention, and you may ask, An expectation or an intention of what? The answer is that it is an expectation or an intention of a proof of that proposition. And Kolmogorov phrased essentially the same explanation by saying that a proposition expresses a problem or task (Ger.\ \emph{Aufgabe}). Soon afterwards, there appeared yet another explanation, namely, the one given by Gentzen, who suggested that the introduction rules for the logical constants ought to be considered as so to say the definitions of the constants in question, that is, as what gives the constants in question their meaning. What I would like to make clear is that these four seemingly different explanations actually all amount to the same, that is, they are not only compatible with each other but they are just different ways of phrasing one and the same explanation. (Martin-L\"of (1987, 410.))
\end{quote}
This passage directly recalls the view expressed in his of \emph{Intuitionistic Type Theory} (1984, 24), that the introduction rules define the meaning:
\begin{quote}
The introduction rules say what are the canonical elements (and equal canonical elements) of the set, thus giving its meaning.
\end{quote}
In a more recent account of his verificationist approach, Martin-L\"of refers to ``Wittgenstein's formula" as the so-called ``verification principle of meaning":
\begin{quote}
It is clear from this what ought to be the general explanation of what a proposition is, namely, that a proposition is defined by stipulating how its proofs, more precisely, canonical or direct proofs, are formed. And, if we take the rules by means of which the canonical proofs are formed to be the introduction rules, I mean, if we call those rules introduction rules as Gentzen did, then his suggestion that the logical constants are defined by their introduction rules is entirely correct, so we may rightly say that a proposition is defined by its introduction rules.

Now what I would like to point out is that this is an explanation which could just as well be identified with the verification principle, provided that it is suitably interpreted. Remember first of all what the verification principle says, namely, that the meaning of a proposition is the method of its verification. The trouble with that principle, considered as a formula, or as a slogan, is that it admits of several differ- ent interpretations, so that there arises the question: how is it to be interpreted? Actually, there are at least three natural interpretations of it. On the first of these, the means of verifying a proposition are simply identified with the introduction rules for it, and there is then nothing objectionable about Wittgenstein's formula, provided that we either, as I just did, replace method by means, which is already plural in form, or else make a change in it from the singular to the plural number: the meaning of a proposition is the methods of its verification. Interpreted in this way, it simply coincides with the intuitionistic explanation of what a proposition is, or, if you prefer, the Gentzen version of it in terms of introduction rules. (Martin-L\"of (2013, 6))
\end{quote}
If meant as `meaning is determined by proof', this does not take into consideration Wittgenstein's self-corrections.\footnote{\S 10, Part VII, p.\ 367 and \S 11, Part VI, p.\ 313. of \emph{Remarks on the Foundations of Mathematics}}

By contrast, in a later phase, Dummett himself admits the possibility of formulating a theory of meaning for logical constants based on assuming that the elimination rules  give meaning to each constant, as noted by P.\ Schroeder-Heister (2018):
\begin{quote}
Most approaches to proof-theoretic semantics consider introduction rules as basic, meaning giving, or self-justifying, whereas the elimination inferences are justified as valid with respect to the given introduction rules. (...)

One might investigate how far one gets by considering elimination rules rather than introduction rules as a basis of proof-theoretic semantics. Some ideas towards a proof-theoretic semantics based on elimination rather than introduction rules have been sketched by Dummett (1991, Ch.\ 13), albeit in a very rudimentary form.
\end{quote}
An earlier statement of the same line of reasoning is found in:
\begin{quote}
Proof-theoretic semantics has several roots, the most specific one being Gentzen's (1934) remarks that the introduction rules in his calculus of natural deduction define the meanings of logical constants while the elimination rules can be obtained as a consequence of this definition. More broadly, it belongs to the tradition according to which the meaning of a term has to be explained by reference to the way it is used in our language. (Schroeder-Heister (2006, 525))
\end{quote}
These remarks have all attempted to 
insist on justifying the so-called `proof-theoretic semantics' by recalling Wittgenstein's `meaning as use' paradigm. In a recent publication, Prawitz recalls Gentzen's dictum when attempting to define what is a valid argument:
\begin{quote}
Gentzen's idea elaborated: To know the \emph{meaning} of a sentence $A$ in $c$-form is to know that \emph{introductions} of $c$ (that is applications of the introduction rule for $c$) are the \emph{canonical} ways of inferring $A$, which is to say: (1) they are valid inferences, and (2) if $A$ can be proved, it can be proved in that way, that is, by a proof whose final step is an application of the introduction rule for $c$ (which we therefore call a \emph{canonical proof}.). (Prawitz (2019, 221))
\end{quote}
On the other hand, Martin-L\"of investigates the application of Gentzen's principle\footnote{From Gentzen (1935) (published in von Plato (2017,153):
\begin{quote}
The ``introductions" present, so to say, the ``definitions" of the signs in question, and the ``eliminations" are actually just consequences thereof, expressed more or less as follows: In the elimination of a sign, the proposition the outermost sign of which is in question, must ``be used only as what it means on the basis of the introduction of this sign."
(...)
I think one could show, by making precise this idea, that the E-inferences are, through certain conditions, \emph{unique} consequences of the respective I-inferences.
\end{quote}
} for assessing the ways in which one can give meaning to Frege's assertion sign `$\vdash$':
\begin{quote}
(...) we may already formulate what it is natural to call the correctness condition for assertion, namely the condition under which it is right, and here several terms are possible to use: right, correct, proper. I am going to use them in the same sense. The condition under which it is right, or correct, or proper, to make an assertion is that you know how to perform the task which constitutes the content of the assertion. This is what I have called the correctness condition for assertion in my abstract.

For acts in general it is usually illuminating to ask, What is the purpose of the act? In this case, if we accept the correctness condition that I just gave, What is the purpose of making an assertion? Then we have already to bring in that the speech act involves not only the speaker, but also the hearer, the receiver of the speech act. So, what is it that the assertor wants to achieve, what is the purpose of making an assertion? Well, if we stick to this knowledge account of assertion that I am discussing right now, then the purpose is nothing but to convey to the hearer that the speaker knows how to fulfil the content, the task which makes up the content. The speech act of assertion has no other purpose than to transmit from the speaker to the hearer the information that the speaker knows how to fulfil the task which makes up the content of the assertion, and this succeeds because the speaker must adhere to the correctness condition for assertion that I just formulated. (Martin-L\"of (2019, 229f))
\end{quote}
Note that here the `Peircean' connection between meaning and purpose via a dichotomy between a speaker and a hearer is rather different from the verificationist perspective advocated earlier. After recalling Gentzen's suggestion in:
\begin{quote}
I began by saying that this whole
lecture will be roughly about what the meaning is of the assertion sign. We are used to the fact that when we ask for the meaning of some linguistic construction, it should be visible somehow from the rules that govern that construction, in general Wittgensteinian terms. The first example of this is of course Gentzen's suggestion that the logical operations are defined by their introduction rules. (Martin-L\"of (2019, 234))
\end{quote}
Martin-L\"of concludes by acknowledging that in this case the elimination rules, rather than the introduction rules, are `meaning-determining for the assertion sign':
\begin{quote}
What about the assertion sign? If you did not have this new rule ($C$-elim), you would only have the usual rules of inference, which are of the form ($C$-intro). If you were to take the assertion sign to be determined by these rules, the assertion $\vdash C$ could not mean anything than that $C$ has been demonstrated, has been inferred by the usual inference rules. And that is not how Frege introduced the assertion sign, what Frege meant by the assertion sign. I explained that earlier on: it is the acknowledgement of the truth of a content that the assertion sign expresses. So, we simply cannot explain the assertion sign by referring to the rules governing it if you only have the rules ($C$-intro). But now we are in a better situation, because we also have the rules ($C$-elim), and they are precisely the rules that are meaning-determining for the assertion sign. (id. ibid.)
\end{quote}
Observe here that Martin-L\"of comes close to Peirce's `\emph{Utterer} vs \emph{Interpreter}' by talking about what he calls a `deontological' approach to meaning via a dichotomy `\emph{Speaker} vs \emph{Hearer}':
\begin{quote}
(...) now I come to the commitment account of assertion, which has its origin in Peirce's work during a very early stage of the last century, 1902-03, I think. Peirce's view was that an assertion should be understood as a taking on of responsibility, taking responsibility for the content of the assertion. (id., 230)
\end{quote}
By highlighting the primary importance of an interaction between a speaker and a hearer, Martin-L\"of acknowledges the role of explaining consequences: what  the hearer can immediately deduce from what the speaker asserts is central to the definition of meaning.

Wittgenstein's clear emphasis on the consequences that can be drawn from a proposition as explaining its meaning does not appear to have received the attention it seems to deserve from advocates of a use-based semantic theory who insist upon the primacy of assertibility conditions. Already in a  1930--31 manuscript Ms-109 V. \emph{Bemerkungen ||} (German; 1930-31 partly published in \emph{Philosophical Remarks}, Foreword), Wittgenstein refers explicitly to his concerns about the view that the meaning of a proposition is given by the way it is verified:
\begin{quote}
Wenn ich sage ,,der Sinn eines Satzes ist dadurch bestimmt, wie er zu verifizieren ist", was mu\ss\ ich dann von dem Sinn des Satzes, da\ss\
dieser Satz || dieses Bild die {\"U}bersetzung || das Portrait jenes Gegenstandes sein soll, sagen? Wie ist das denn zu verifizieren?\footnote{which may be translated as:
\begin{quote}
If I say ``the sense of a proposition is determined by how it is to be verified," what must I think of the sense of the proposition that
this sentence || this picture the translation || the portrait of that object is supposed to be? How is that to be verified?
\end{quote}
Ms-109,77,
Wittgenstein \emph{Nachlass} Ms-109: V, \emph{Bemerkungen} (WL) -  (IDP).
}
\end{quote}
The need for an explanation, rather than an interpretation, appears rather often in an earlier (1930--31) manuscript:
\begin{quote}
Ich will also eigentlich sagen: Es gibt nicht Grammatik \& Interpretation der Zeichen. Sondern soweit von einer Interpretation, also von einer Erkl\"{a}rung der Zeichen, die Rede sein kann, soweit mu\ss\ sie die Grammatik selbst besorgen.

     Denn ich brauchte nur zu fragen: Soll die Interpretation durch S{\"a}tze erfolgen? Und in welchem Verh\"{a}ltnis sollen diese S\"{a}tze zu der Sprache stehen die sie schaffen?\footnote{which may be translated as:
 \begin{quote}
 What I really want to say is: There is no such thing as grammar and interpretation of the signs. Rather, so far as there can be talk of an interpretation, that is, of an explanation of the signs, the grammar itself must provide it.
 
     Because I only had to ask: Should the interpretation be done by propositions? And how should these sentences relate to the language they create?
 \end{quote}
Ms-109,129 (Wittgenstein \emph{Nachlass} Ms-109: V, Bemerkungen (WL) - (IDP).)}
\end{quote}

 Neither has there been much concern with this particular aspect of Wittgenstein's very late writings in most of the literature on Wittgenstein's philosophy, to the best of our knowledge. Thus, to develop a reformulation of the so-called Dummett--Prawitz theory of meaning for the language of mathematics, which was called by both authors a `Wittgensteinian theory', it seems advisable to bring forward the actual perspective of Wittgenstein on language, meaning and use.

\section{Later writings}

Already in a manuscript written in 1937--38, one can find remarks insisting on the connection between the meaning and the use of a word, and the consequences one can draw from a sentence where the word is used:
\begin{quote}
``Diese Vorstellung ist doch `rot'\ " als k\"onnte ich blo\ss\ durch einen Akt w\"ahrend der Definition von `rot' unterscheiden ob ich die Vorstellung oder den K{\"o}rper, oder seine Farbe, mit ``rot" meine. Freilich gibt es ein: die Farbe eintrinken und anderes; aber wie hilft mir das beim Definieren? Inwiefern habe ich der Farbe den Namen gegeben, wenn ich auf die Farbe beim Definieren meine
Aufmerksamkeit gerichtet hatte? Die Definition ist doch zu k{\"u}nftigem Gebrauch da, ? \& wie n\"utzt dabei jener seelische Akt? Man k{\"o}nnte sagen: Was Du beim Definieren meinst, wenn dies Meinen ein seelischer Akt sein soll, ist f\"ur den Gebrauch \& daher f\"ur die Bedeutung der Definition ganz gleichg{\"u}ltig. Es wird drauf ankommen, was Du mit dieser Definition machst. Was Du also ``rot" benannt hast, zeigt der Kalk{\"u}l mit diesem Wort. Ich verstehe nat{\"u}rlich, da\ss\ die verschiedenen Arten des Zeigens \& des Anwendens der Aufmerksamkeit beim Definieren n{\"u}tzlich sein kann; z.B. zum Einpr{\"a}gen der Definition. Aber es ist ein Aberglaube, da\ss\ man durch einen bestimmten seelischen Akt beim
Definieren nicht z.B.\ der Wange sondern der Empfindung den Namen ``Zahnschmerz" gibt. Als f\"uhrte man den Namen dadurch diesem \& nicht jenem Gegenstande zu. Als erhielte ihn nun der eine \& nicht der andere. Als k{\"o}nne man beim Definieren durch einen Akt des `Meinens' diesem \& nicht jenem den Namen anheften. --  (Eine Fragestellung scheint hier falsch) Damit ist nicht zu verwechseln: ich kann durch einen seelischen Akt eine || diese assoziative Verbindung herstellen \& nicht eine andere. Aber: ob es mir gelungen ist die \& die Verbindung herzustellen, zeigen die Folgen. D.h. ob ich mir etwas durch diese Art des Ansehens einpr{\"a}ge zeigt sich, wenn mein Ged{\"a}chtnis auf die Probe gestellt wird.\footnote{which may be translated as:\begin{quote}
``But this idea is `red'\ " as if by just one act during the definition of `red' I could distinguish whether by `red' I mean the idea or the body, or its colour. Of course there is one: drinking the paint and other such; but how does that help me define? To what extent did I name the color when I mean the color is defining? The definition is there for future use, -- \& what use is that mental act? One could say: What you mean when defining, if this meaning is supposed to be an act of the soul, is completely irrelevant for the use \& therefore for the meaning of the definition. It will depend on what you do with that definition. So what you named ``red" shows the calculus with this word. I understand of course that the different ways of showing \& applying attention can be useful in defining; e.g.\ to memorize the definition. But it is a superstition that through a certain mental act one can
not define e.g.\ the cheek but give the sensation the name ``toothache". As if the name was thereby attributed to this and not that object. As if the one and not the other received it. As if one could attach the name to this \& not that when defining by an act of `intending'. ? This is not to be confused with: I can, through an act of the soul, create a || make this associative connection \& not another. But: whether I managed to establish such \& so connection is shown by the consequences. I.e.\ whether I memorize something is shown when my memory
is put to the test.
\end{quote}
(Ms 120,2v,3r,3v,4r)
Wittgenstein \emph{Nachlass} Ms-120: XVI (WL) -  (IDP).}

\end{quote}

Notes from a series of \emph{Lectures on Religious Belief}, given in Cambridge sometime around summer 1938, taken by Y.\ Smythies and published in \emph{Lectures and Conversations on Aesthetics, Psychology and Religious Belief}, one finds passages linking use and consequences:

\begin{quote}
If you say to me -- ``Do you cease to exist?'' -- I should be bewildered, and
would not know what exactly this is to mean. ``If you don't cease to exist,
you will suffer after death'', there I begin to attach ideas, perhaps ethical
ideas of responsibility. The point is, that although these are well-known words,
and although I can go from one sentence to another sentence, or to pictures,
[{I don't know what consequences you draw from this statement\/}].\footnote{p.\ 69--70 of the {\it Lectures and Conversations on Aesthetics, Psychology and Religious Belief\/}.}

``Yes, this might be a disagreement---{if he himself were
to use the word in a way in which I did not expect, or were to draw conclusions
I did not expect him to draw\/}.''\footnote{Ibid., 71.}

When I say he's using a picture I'm merely making a {\em grammatical\/} remark.
{[What I say] can only be verified by the consequences he does or does not
draw.}\footnote{Ibid., 72.}
\end{quote}

\begin{quote}
``... someone who does philosophy or psychology will perhaps say ``{\it I\/}
feel that I think in my head''. But what that means he won't be able to say. For
he will not be able to say {\it what\/} kind of feeling that is; but merely to
use the expression that he `feels'; as if he were saying ``{\it I\/} feel this
stitch {\it here\/}''. Thus he is unaware that {it remains to be
investigated what his expression ``I feel'' means here, that is to say:
what consequences [{\it Konsequenzen\/}] we are permitted
to draw from this utterance.}
Whether we may draw the same ones as we would from the utterance ``I feel a
stitch here''.\footnote{\S 350, p.\ 69$^{\rm e}$ of the {\it Remarks on the Philosophy
of Psychology\/}, Vol. I.}
\end{quote}

\begin{quote}
For the question is not, `What am I doing when . . .?' (for
this could only be a psychological question) -- but rather, `What meaning does
the statement have, what can be deduced from it, {what consequences does it
have\/}?'\footnote{\S 38, p.\ 8$^{\rm e}$ of the {\it Remarks on
the Philosophy of Psychology\/}, Vol. II. (TS 232, bearing remarks
which stem from manuscripts written between November 1947 and August 25th 1948)}
\end{quote}

\begin{quote}
I {\sc know} that this is my foot. I could not accept any
experience as proof to the contrary.---That may be an exclamation;
{but what {\it follows\/} from it?} At least that I shall act with a
certainty that knows no doubt, in accord with my belief.\footnote{\S 360, p.\ 47$^{\rm e}$ of {\it On Certainty/ \"Uber Gewissheit\/}.}
\end{quote}

\begin{quote}
One is often bewitched by a word. For example, by the word
``know''.

Is God bound by our knowledge? Are a lot of our statements {\it incapable\/} of
falsehood? For that is what we want to say.

I am inclined to say: ``That {\it cannot\/} be false.'' That is interesting;
{but what consequences has it?}\footnote{Ibid., \S \S 435--7, p.\ 57$^{\rm e}$. (dated 27.3.51)}
\end{quote}

In the first volume of the \emph{Last Writings on the Philosophy of Psychology}, which is not based on a typescript as is the case for the previously published \emph{Remarks on the Philosophy of Psychology}, but on manuscript writings dating from the period between 22 October 1948 and 22 March 1949, except for the last remark dated 20 May (second half of MS 137 and the whole of MS 138), one can find examples of very explicit references to the connections between the use of words and the consequences of the corresponding utterance, sometimes presented positively, as in the following:
\begin{quote}
What are you telling me when you {\it use\/} the words . . .?
What can I do with this utterance? {What consequences does it
have?\/}\footnote{\S 624, p.\ 80$^{\rm e}$ of {\it Last Writings on the Philosophy of Psychology I\/}.}
\end{quote}

\begin{quote}
What does anyone tell me by saying ``Now I see it as ...''?
{What consequences has this information?\/} What can I do with
it?\footnote{Id.\ Ibid., \S 630.  (which also appears in the {\it Investigations\/}, Part II,
section xi, p.\ 202$^{\rm e}$)}
\end{quote}

\begin{quote}
The report ``The word . . . was crammed full of its meaning''
{\it is used quite differently\/}, {has quite different consequences\/},
from ``It had the meaning ...''.\footnote{Ibid., \S 785, p.\ 100$^{\rm e}$.}
\end{quote}

The renewed examination of Wittgenstein's \emph{Nachlass} presented here corroborates and augments the significance of Wittgenstein's suggestion that meaning is determined by explicating immediate consequences of a statement or term. These new findings counter the claim that the so-called Gentzen-Dummett-Prawitz-Martin-L\"of `meaning as determined by assertability conditions / introduction rules' are properly designated as `Wittgensteinian views'.
Instead, these passages strongly suggest connections with the ``pragmatist"/``dialogical" approaches to meaning, already in Peirce's writings on the interaction between the \emph{Interpreter} and the \emph{Utterer}, which likewise appear to bear upon the ``game"/``dialogue" approaches to meaning (Lorenzen, Hintikka, Fra\"{\i}ss\'e). The \emph{Inversion Principle}\footnote{To recall Prawitz (1965):
\begin{quote}
{\sc Inversion Principle}. Observe that the elimination rule is, in a sense, the inverse of the corresponding introduction rule: by an application of an elimination rule one essentially restores what had already been established if the major premiss of the application was inferred by an application of an introduction rule. This relationship between the introduction rules and the elimination rules is roughly expressed by the following principle, which I shall call \emph{inversion principle}:

\emph{Let $\alpha$ be an application of an elimination rule that has $B$ as consequence. Then, deductions that satisfy the sufficient condition (in the list above) for deriving the major premiss of $\alpha$ (if any), already ``contain" a deduction of $B$; the deduction of $B$ is thus obtainable directly from the given deductions without the addition of $\alpha$.} (Prawitz (1965, 33))
\end{quote}
In a footnote Prawitz refers to Lorenzen's principle of inversion:
\begin{quote}
What Lorenzen has called the principle of inversion (see e.g.\ Lorenzen (1950) and (1955))) is closely related to this idea emanating from Gentzen. (For a correct statement of Lorenzen's principle, see Hermes (1959).) (Id.\ Ibid.)
\end{quote}
} 
 (advocated by Gentzen and Prawitz) also appears to have the same kind of ``metalevel" explication of meaning as those game/dialogical approaches, including Tarski's ``metalanguage" approach. Hence in Natural Deduction, the reduction rules (which formalise the \emph{Inversion Principle}) should be seen as meaning-giving, not a mere justification of the elimination rule(s) based on the introduction rule(s). The reduction rules of labelled natural deduction (`\emph{Grundgesetze} alongside \emph{Begriffsschrift}' (de Queiroz (1992), de Queiroz and Gabbay (1999), de Queiroz et al.\ (2011))) examined below, extends this approach to `Identity types' which correspond to propositional equality. Significantly, the present examination of these materials make an important step towards a formal counterpart to the ``meaning is use" dictum, while highlighting an important common thread from Wittgenstein's very early to very late writings.

\section{Lorenzen's Dialogue Semantics and Hintikka's Games Semantics}
There are at least two research programs dealing with the interface of meaning, knowledge, and logic in the context of dialogues, games, or more generally interaction. One of these is also an alternative perspective on proof theory and meaning theory, advocating that Wittgenstein's ``meaning as use" paradigm, as understood in the context of proof theory, by which the so-called reduction rules (showing the workings of elimination rules on the result of introduction rules) should be seen as appropriate to formalise the explanation of the (immediate) consequences one can draw from a proposition, thus to show the function/purpose/usefulness of its main connective in the calculus of language. To recall D.\ Prawitz' \emph{Natural Deduction} (1965) original formulation: `reduction steps'  are defined as rules which operate on proofs, thus they constitute `metalevel' rules rather than `object level' rules. This suggests they can indeed be seen as playing a `semantical' role, in analogy to several other definitions (Tarski's truth conditions, Lorenzen's dialogical rules, Hintikka's game semantics rules, Abelard-Elo\"{\i}se evaluation game rules\footnote{Cf.\ V\"an\"a\"anen (2022):
\begin{quote}
The game $G(M, \phi)$
has two players Abelard and Eloise. Intuitively, Eloise defends the proposition that
$\phi$ is (informally) true in $M$ and Abelard doubts it. All through the game the players
are inspecting an $L$-sentence and an assignment $s$. (...)

Suppose the position is $(\psi, s)$.\\
(1) If $\psi$ is a literal, the game ends and Eloise wins if $\psi\in \Gamma(s)$. Otherwise
Abelard wins.\\
(2) If $\psi$ is $\psi_0\land\psi_1$, then Abelard chooses whether the next position is $(\psi_0, s)$ or
$(\psi_1, s)$.
(3) If $\psi$ is $\psi_0 \lor\psi_1$, then Eloise chooses whether the next position is $(\psi_0, s)$ or
$(\psi_1, s)$.\\
(4) If $\psi$ is $\forall x\theta$, then Abelard chooses $a \in M$ and the next position is $(\theta, s(a/x))$.\\
(5) If $\psi$ is $\exists x\theta$, then Eloise chooses $a \in M$ and the next position is $(\theta, s(a/x))$.\\
(...) The first to formulate the Evaluation Game explicitly was probably J.\ Hintikka [(1968)]
 who later advocated the importance and usefulness of the game forcefully. Hintikka mentions L.\ Wittgenstein [\emph{Investigations}] as an inspiration. Earlier L.\ Henkin [(1961)] formulated a game theoretic approach to quantifiers pointing out generalizations to very
general infinitary and partially ordered quantifiers. Nowadays the Evaluation Game
is a standard tool in mathematical, philosophical and computer science logic. 
\end{quote}
}).

Regarding the interface of meaning, logic and interaction, E.\ Krabbe (2001) allows for a perspective under which Lorenzen's dialogical acts could be seen as an implementation of Wittgenstein's meaning-is-use for the language of mathematics:
\begin{quote}
Thus it appears that there are two different objectives for dialogue theory in the context of foundations of mathematics. The first and primary objective is to provide a criterion of meaningfulness in terms of dialogical acts that are often described as attacks and defenses. The well-known slogan `Don't ask for the meaning, ask for the use' is here translated as `To know what a sentence means, you must know how to attack and defend it'.
\end{quote}
In his own contribution to the same volume, `A  Sceptical Look', Wilfrid Hodges (2001) starts by quoting C.\ S.\ Peirce's early origin of the duality between the game rule for the universal quantifier and that for the existential quantifier:
\begin{quote}
In his second Lowell Lecture in 1898, C.S.\ Peirce explains the difference between `every' and `some':
\begin{quote}
When I say `every man dies', I say you may pick out your man for yourself and provided he belongs to \emph{this here} world you will find he will die. The `some' supposes a selection from `this here' world to be made by the \emph{deliverer} of the proposition, or made in his interest. The `every' \emph{transfers} the function of selection to the \emph{interpreter} of the proposition, or to anybody acting in his interest.
\end{quote}
This is interesting and it deserves some unpicking. Peirce imagines two people, a deliverer and an interpreter. (...) At first sight it seems that the parallel is unconvincing; in the case of `every', the interpreter can pick any man, but in the case of `some' the deliverer must pick a suitable man. But if we view it that way, it's irrelevant who makes the choice, and Peirce never had any reason to distinguish the deliverer from the interpreter. As often, Peirce leaves it to us to decide what point he was trying to make.
\end{quote}
Notice that the often named `pragmatist' meaning theory (as sometimes assigned to Dummett's later views) advocates that the elimination rules of the logical connective, which, in some (incomplete) sense, expresses the immediate consequences one can draw from the sentence, finds its origin exactly in the best known pragmatist, namely, C.\ S.\ Peirce, and this is no coincidence. It should be no surprise to see that rules for dialogical logic recalls Peirce's perspective:
\begin{quote}
A particle rule (also known as an argumentation form) abstractly describes the way a formula of a given main connective may be objected to, and how to answer the objection. By definition, an argumentation form is a tuple consisting of (1) an assertion, (2) a set of attacks, (3) a set of defenses, and (4) a relation specifying for each attack the corresponding defense(s). Argumentation forms are abstract in the sense that, in their definition, no reference is made to the context of argumentation in which the rule is applied (understood as a situation in a process of argumentation, e.g., as the set of argument moves made before the assertion for which the particle rule is defined). The particle rules thus constitute the local semantics of a logic, for they determine the dialogical meaning of each logical constant but say nothing about the way this meaning may be related to anything else. Yet the local semantics given by the particle rules are intrinsically dependent on the notion of dialogue, since they describe sequences of language acts. (Keiff (2009)).
\end{quote}

Indeed, a number of authors have drawn attention to this origin of game-based semantics such as Lorenzen's and Hintikka's. For the latter, Pietarinen (2014) acknowledges:
\begin{quote}
Jaakko Hintikka is the undisputed modern achitect of \emph{game-theoretical semantics} (GTS), but the spirit and often even the letter of that theory date back to polymath Charles Peirce and his groundbreaking discoveries in the logical studies of the 19$^{\mathrm{th}}$ and the early 20$^{\mathrm{th}}$-century.\end{quote}
In fact, an alternative to the ``proof-theoretic semantics" (Gentzen, Prawitz, Dummett, Martin-L\"of) which advocates that ``the meaning of a logical connective is given by its assertability conditions" draws on Wittgenstein's early suggestion that the explanation of the consequences which can be drawn from a proposition would be more appropriate for the so-called `meaning is use' paradigm; this appears as early as 1912 in a letter to Russell. Now, when looking at the nature of more recent approaches which can be classified as ``dialogue-like" explanation of meaning, such as Hintikka's game-theoretic semantics and Lorenzen dialogue-games, one can see a parallel to Peirce's ``\emph{Utterer} vs \emph{Interpreter}" framework in the distinction between universal and existential quantifiers. Peirce was concerned with empirical statements, but the dialogical stand suggested by the dichotomy ``\emph{Utterer} vs \emph{Interpreter}" can equally be logical explanations of meaning via explanation of consequences which can be drawn from a proposition.\footnote{Peirce was likewise concerned with semiotic quite generally, and with both logic and mathematics. Specifying terms or concepts,  expect any fundamental difference between empirical and logical or mathematical specifications of terms. Cf.\ Peirce's remark:
\begin{quote}
In studying logic, you hope to correct your present ideas of what reasoning is good, what bad. This, of course, must be done by reasoning. Some writers fancy that they see some absurdity in this. They say it would be a `\emph{petitio principii}'. Let us rather state the case thus. At present you are in possession of a \emph{logica utens} which seems to be unsatisfactory. The question is whether, using this unsatisfactory \emph{logica utens}, you can make out wherein it must be modiified, and can attain to a better system. This is a truer way of stating the question; and so stated, it appears to pres­ent no such insuperable difficulty, as is pretended. (\emph{Collected Papers}, 2:191).
\end{quote}
\  }

Previously we have pointed out parallels between our proposal and other approaches to the semantics of logical connectives based on explaining the (immediate) consequences of the corresponding proposition, such as Lorenzen's  (1950, 1955, 1969, p.\ 25) dialogical games and Hintikka's (1968, 1979, p.\ 3) semantical games. The general underlying principle is the logical \emph{Inversion Principle}, uncovered by Gentzen, and later by Lorenzen: the elimination procedure is the exact inverse of the introduction, therefore all that can be asked from an assertion is what is indicated  explaining the elimination procedure. In the spirit of Wittgenstein's remark in his (1930--31) Ms 110,61:
 `Die Erkl\"arung ist \"aquivalent mit der Bedeutung' [`The explanation is equivalent to the meaning'].\footnote{Wittgenstein \emph{Nachlass} Ms-110: VI, Philosophische Bemerkungen (WL), (IDP) (accessed 14 May 2022)}
One can look at Lorenzen's dialogical games and try to show how the  \emph{Inversion Principle} resides in the game-approach by comparing it to the reduction-based account of meaning. The correspondence works as follows: (i) the assertion corresponds to the introduction rule(s) (and the respective constructors);
(ii) the attacks correspond to the elimination rule(s) (and the respective destructors);
(iii) the set of defenses correspond to the reduction rules (the effect of elimination rules on the result of introduction rules);
(iv) the relation specifying for each attack the corresponding defense(s) are defined by the result of reduction rules.

Taking from the rules of \emph{reduction} between proofs in a system of Natural Deduction enriched with terms alongside formulas, such as in the so-called Curry-Howard interpretation (Howard (1980)) of which Martin-L\"of's type theory can be seen as an instance,\footnote{Howard (2014) remarks on the so-called {\em Curry-Howard `Formulae-as-Types'}:
\begin{quote}
{[}de Bruijn] discovered the idea of {\bf derivations as terms}, and the accompanying idea of formulae-as-types, on his own.  (...) 

Martin-L\"of's work originated from mine. He has always given me credit and we are good friends.

On further thought, I need to mention that, in that first conversation, Martin-L\"of suggested that the derivations-as-terms idea would work particularly well in connection with Prawitz's theory of natural deduction. I thought: okay, but no big deal. Actually, at that time, I was not familiar with Prawitz's results (or, if at all, then only vaguely). But it was a bigger deal than I had thought, because Prawitz's reductions steps for a deduction correspond directly to reduction steps for the associated lambda term!
\end{quote}
} and isolating the terms corresponding the derivations, one can draw the following parallels:\footnote{
The rules of reduction between proofs in the system of Natural Deduction presented in the style of Curry-Howard `derivation-as-terms' are taken from (de Queiroz et al.\ 2011).}
Looking at the conclusion of reduction inference rules, one can take the destructor as being the Attack (or `Nature', its counterpart in the terminology of Hintikka's Game-Theoretical Semantics), and the constructor as being the Assertion (or Hintikka's `Myself'). This way the game-theoretic explanations of logical connectives find direct counterpart in the functional interpretation with the semantics of convertibility:

\medskip

\noindent $\land$-$\beta$-{\it reduction\/}
$$\displaystyle{{\displaystyle{{a:A \qquad b:B} \over
{\langle a,b\rangle:A\land B}}\land\mbox{\it -intr}} \over
{{\tt FST}(\langle a,b\rangle):A}}\land\mbox{\it -elim} \qquad \qquad
\twoheadrightarrow_\beta \qquad \qquad a:A$$
$$\displaystyle{{\displaystyle{{a:A \qquad b:B} \over
{\langle a,b\rangle:A\land B}}\land\mbox{\it -intr}} \over
{{\tt SND}(\langle a,b\rangle):B}}\land\mbox{\it -elim} \qquad \qquad
\twoheadrightarrow_\beta \qquad \qquad b:B$$
Associated rewritings:\\
${\tt FST}(\langle a,b\rangle)=_\beta a$\\
${\tt SND}(\langle a,b\rangle)=_\beta b$\\
\ 
\\
\begin{tabbing}
\underline{Assertion}/{\it Introd.\/} $\quad$ \= \underline{Attack}/{\it Elim.\/} $\qquad\qquad\quad \twoheadrightarrow_\beta \qquad\qquad \qquad
$ \= \underline{Defense} \\
\\
Conjunction (`$\land$'): \\
${A}\land{B}$ \> $\quad \quad$ $L$? \> ${A}$ \\
${A}\land{B}$ \> $\quad \quad$ $R$? \> ${B}$ \\
\\
$\langle a,b\rangle:{A}\land{B}$ \> {\tt FST}$(\langle a,b\rangle)$ $\qquad \twoheadrightarrow_\beta$ \> $a:{A}$ \\
$\langle a,b\rangle:{A}\land{B}$ \> {\tt SND}$(\langle a,b\rangle)$ $\qquad \twoheadrightarrow_\beta$ \> $b:{B}$ \\
\end{tabbing}

\noindent $\lor$-$\beta$-{\it reduction\/}
$$\displaystyle{{\displaystyle{{a:A} \over {{\tt inl}(a):A\lor B}}\lor\mbox{\it -intr} \  \displaystyle{{[x:A]} \atop {f(x):C}} \  \displaystyle{{[y:B]} \atop {g(y):C}}} \over
{{\tt CASE}({\tt inl}(a),\upsilon x.f(x),\upsilon y.g(y)):C}}\lor\mbox{\it -elim} \ \ \twoheadrightarrow_\beta \  \displaystyle{{a:A} \atop {f(a/x):C}}$$
$$\displaystyle{{\displaystyle{{b:B} \over {{\tt inr}(b):A\lor B}}\lor\mbox{\small\it -intr} \  \displaystyle{{[x:A]} \atop {f(x):C}} \  \displaystyle{{[y:B]} \atop {g(y):C}}} \over
{{\tt CASE}({\tt inr}(b),\upsilon x.f(x),\upsilon y.g(y)):C}}\lor\mbox{\small\it -elim} \ \ \twoheadrightarrow_\beta \ \  \displaystyle{{b:B} \atop {g(b/y):C}}$$
Associated rewritings:\\
${\tt CASE}({\tt inl}(a),\upsilon x.f(x),\upsilon y.g(y))=_\beta f(a/x)$\\
${\tt CASE}({\tt inr}(b),\upsilon x.f(x),\upsilon y.g(y))=_\beta g(b/y)$\\
\ 
\\
\begin{tabbing}
\underline{Assertion}/{\it Introd.\/} $\quad$ \= \underline{Attack}/{\it Elim.\/} $\qquad\qquad\quad \twoheadrightarrow_\beta \qquad\qquad \qquad
$ \= \underline{Defense} \\
\\

Disjunction (`$\lor$'): \\
${A}\lor{B}$ \> $\quad \quad$ ? \> ${A}$ \\
${A}\lor{B}$ \> $\quad \quad$ ? \> ${B}$ \\
\\
{\tt inl}$(a):{A}\lor{B}$ \> {\tt CASE}({\tt inl}$(a),\upsilon x.f(x),\upsilon y.g(y))$ $\quad \twoheadrightarrow_\beta$ \> $a:A$, $f(a/x):{C}$ \\
{\tt inr}$(b):{A}\lor{B}$ \> {\tt CASE}({\tt inr}$(b),\upsilon x.f(x),\upsilon y.g(y))$ $\quad \twoheadrightarrow_\beta$ \> $b:B$, $g(b/y):{C}$ \\
\\
\end{tabbing}
Whilst in both {\tt FST}$(\langle a,b\rangle)$ and {\tt SND}$(\langle a,b\rangle)$ the destructors allow access to either of the conjuncts, in {\tt CASE}({\tt inl}$(a),\upsilon x.f(x),\upsilon y.g(y))$ and
 {\tt CASE}({\tt inr}$(b),\upsilon x.f(x),\upsilon y.g(y))$ the destructor is not given access to the either disjunct but must ask for whichever disjunct comes from the introduction.

\medskip

\noindent $\rightarrow$-$\beta$-{\it reduction\/}
$$\displaystyle{{\displaystyle{\ \atop {a:A}} \qquad \displaystyle{{\displaystyle{{[x:A]} \atop {b(x):B}}} \over
{\lambda x.b(x):A\rightarrow B}}\rightarrow\mbox{\it -intr}} \over
{{\tt APP}(\lambda x.b(x),a):B}}\rightarrow\mbox{\it -elim} \quad \twoheadrightarrow_\beta
\quad \displaystyle{{a:A} \atop {b(a/x):B}}$$
Associated rewriting:\\
${\tt APP}(\lambda x.b(x),a)=_\beta b(a/x)$\\
\ 
\\
\begin{tabbing}
\underline{Assertion}/{\it Introd.\/} $\quad$ \= \underline{Attack}/{\it Elim.\/} $\qquad\qquad\quad \twoheadrightarrow_\beta \qquad\qquad \qquad
$ \= \underline{Defense} \\
\\
Implication (`$\rightarrow$'): \\
${A}\rightarrow{B}$ \> $A$ $\quad$ ? \> {$B$} \\
\\
$a:A$\\
$\lambda x.b(x):A\to B$ \> {\tt APP}$(\lambda x.b(x),a)$ $\qquad
\twoheadrightarrow_\beta$ \> $b(a/x):{B}$ \\
\\
\end{tabbing}

\noindent $\forall$-$\beta$-{\it reduction\/}
$$\displaystyle{{\displaystyle{\ \atop {a:D}} \qquad \displaystyle{{\displaystyle{{[x:D]} \atop {f(x):P(x)}}} \over
{\Lambda x.f(x):\forall x^D.P(x)}}\forall\mbox{\it -intr}} \over
{{\tt EXTR}(\Lambda x.f(x),a):P(a)}}\forall\mbox{\it -elim} \quad \twoheadrightarrow_\beta
\quad \displaystyle{{a:D} \atop {f(a/x):P(a)}}$$
Associated rewriting:\\
${\tt EXTR}(\Lambda x.f(x),a)=_\beta f(a/x)$\\
\

\noindent $\exists$-$\beta$-{\it reduction\/}
$$\displaystyle{{\displaystyle{{a:D \quad f(a):P(a)} \over
{\varepsilon x.(f(x),a):\exists x^D.P(x)}}\exists\mbox{\it -intr} \quad
\displaystyle{{[t:D,g(t):P(t)]} \atop {d(g,t):C}}} \over
{{\tt INST}(\varepsilon x.(f(x),a),\sigma g.\sigma t.d(g,t)):C}}\exists\mbox{\it -elim} \  \twoheadrightarrow_\beta\ 
\displaystyle{{a:D,f(a):P(a)} \atop {d(f/g,a/t):C}}$$
Associated rewriting:\\
${\tt INST}(\varepsilon x.(f(x),a),\sigma g.\sigma t.d(g,t))=_\beta d(f/g,a/t)$\\
\

\begin{tabbing}
\underline{Assertion}/{\it Introd.\/} $\quad$ \= \underline{Attack}/{\it Elim.\/} $\qquad\qquad\quad \twoheadrightarrow_\beta \qquad\qquad \qquad
$ \= \underline{Defense} \\
\\
Universal Quantifier (`$\forall$'): \\
$\forall x^{D}.{P}(x)$ \> $a:{D}$ ? \> ${P}(a)$ \\
\\
$a:A$\\
$\Lambda x.f(x):\forall x^{D}.{P}(x)$ \\
\> {\tt EXTR}$(\Lambda x.f(x),a)$ $\qquad
\twoheadrightarrow_\beta$ \> $f(a/x):{P}(a)$ \\
\\
Existential Quantifier (`$\exists$'): \\
$\exists x^{D}.{P}(x)$ \> $\quad \quad$ ? \> $a:{D}, {P}(a)$ \\
\\
\> \> $a:D, f(a):P(a)$\\
$\varepsilon x.(f(x),a):\exists x^{D}.{P}(x)$ \\
\> {\tt INST}$(\varepsilon x.(f(x),a),\sigma g.\sigma t.d(g(t),t))$ $\ \ \twoheadrightarrow_\beta$ \> $d(f(a)/g(t),a/t):{C}$ \\
\\
\end{tabbing}
Whilst in {\tt EXTR}$(\Lambda x.f(x),a)$  the destructor has access to the witness (can use a generic element)
$a$, in {\tt INST}$(\varepsilon x.(f(x),a),\sigma g.\sigma t.d(g(t),t))$ the only option is to eliminate over the witness $a$ which is inside the term $\varepsilon x.(f (x), a)$ built with the constructor because it was chosen by the introduction.

\medskip

\noindent $Id$-$\beta$-{\it reduction\/}
$$\displaystyle{{\displaystyle{{a=_r b:A} \over {r(a,b):Id_A(a,b)}}Id\mbox{\it -intr\/} \qquad
\displaystyle{{[a=_t b:A]} \atop {d(t):C}}} \over
{{\tt REWR}(r(a,b),\sigma t.d(t)):C}}Id\mbox{\it -elim\/} \qquad
\twoheadrightarrow_\beta \qquad
\displaystyle{{a=_r b:A} \atop {d(r/t):C}}$$
Associated rewriting:\\
${\tt REWR}(r(a,b),\sigma t.d(t))=_\beta d(r/t)$\\
\
\begin{tabbing}
\underline{Assertion}/{\it Introd.\/} $\quad$ \= \underline{Attack}/{\it Elim.\/} $\qquad\qquad\quad \twoheadrightarrow_\beta \qquad\qquad \qquad
$ \= \underline{Defense} \\
\\
Propositional Equality (`$Id_A(a,b)$'): \\
$Id_A(a,b)$ \> $\quad \quad$ ? \> $a=_r b:A$ \\
\\
$r(a,b);Id_A(a,b)$ \> {\tt REWR}$(r(a,b),\sigma t.d(t))$ $\ \ \twoheadrightarrow_\beta$ \> $d(r/t):{C}$ \\

\end{tabbing}

In {\tt REWR}$(r(a,b),\sigma t.d(t))$ the destructor has no choice regarding the 
`reason' for $a$ being equal to $b$, since $r$ will have been chosen by the time of the application of the introduction rule.

\medskip

Drawing these parallels is intended to exhibit a connection with the ``pragmatist"/``dialogical" approaches to meaning as revealed in Peirce's writings on the interaction between the \emph{Interpreter} and the \emph{Utterer} which seem to pertain to the ``game"/``dialogue" approaches to meaning (Lorenzen, Hintikka, Ehrenfeucht-Fra\"{\i}ss\'e\footnote{Of course, the main contribution of Ehrenfeucht-Fra\"{\i}ss\'e methods concern their use of the notion of \emph{winning strategy} in order to prove/disprove elementary equivalence between two models, but the actual rules of the games between \emph{Player I} and \emph{Player II} were based on the distinction between universal and existential formulas in the sense of which player is to choose the next element from a structure. Such a distinction recalls Peirce's \emph{Utterer/Interpreter}, Lorenzen's \emph{Attack/Defense}, Hintikka's \emph{Nature/Myself}.}). 

\section{Conclusions}
Using Wittgenstein's \emph{Nachlass} (as we had not done previously) and examining those passages quoted here from his unpublished works, strongly corroborate and augment the main thrust of our continuing research on the significance of Wittgenstein's suggestion that meaning is specified by explaining the immediate consequences of a term or statement. 
Examining these new findings in the \emph{Nachlass}, we have sought to connect Wittgenstein's account of 'meaning as use', specified by the 'calculus' involved in specifying or using terms (or statements), with the pragmatist approaches to meaning as revealed in Peirce's writings on the interaction between the \emph{Interpreter} and the \emph{Utterer} which seem to bear on the ``game"/``dialogue" approaches to meaning (Lorenzen, Hintikka, Fra\"{\i}ss\'e).

In this connection, we have also drawn attention to the fact that the so-called \emph{Inversion Principle} (as put forward by Gentzen and Prawitz)  has the same status as ``metalevel" explanation of meaning as do these game/dialogical approaches, including Tarski's ``metalanguage" approach. 
So, in Natural Deduction, the reduction rules (which formalise the \emph{Inversion Principle}) should be seen as meaning-giving, not merely as a justification of the elimination rule based on the introduction rule.
The reduction rules of labelled natural deduction include the case of the ``Identity types" which were crucial for the more mathematically based recent publications connecting proof theory and homotopy theory.

Most importantly, we believe this material makes an important step towards a formal counterpart to the ``meaning is use" dictum, despite several slight variations on this theme
there is a common thread from the early writings (\emph{Notebooks}, \emph{Letters to Russell} (1912)) to the texts of the so-called second Wittgenstein, whether transitional or very late. By way of further illustrating the main point, now quoting from Wittgenstein's own writings:
\begin{quote}
Was hei{\ss}t es: den Goldbach'schen Satz glauben? Worin besteht dieser Glaube? In einem Gef{\"u}hl der Sicherheit, wenn wir den Satz aussprechen, oder h{\"o}ren? Das interessiert uns nicht. Ich wei\ss\ ja auch nicht, wie weit dieses Gef{\"u}hl durch den Satz selbst hervorgerufen sein mag. Wie greift der Glaube in diesen Satz ein? Sehen wir nach, welche Konsequenzen er hat, wozu er uns bringt. ``Er bringt mich zum Suchen nach einem Beweis dieses Satzes". -- Gut, jetzt sehen wir noch nach, worin Dein Suchen eigentlich besteht; dann werden wir wissen, ``wie es sich mit Deinem Glauben an den Satz verh{\"a}lt. ||  was es mit dem Glauben an den Satz auf sich hat.\footnote{which may be translated as
\begin{quote}
``What does it mean: to believe Goldbach's theorem? What is this belief in? In a feeling of security when we say or hear the sentence? That's of no interest here. I don't know how far this feeling may have been evoked by the sentence itself. How does belief intervene in this sentence? Let's see what consequences it has, what it leads us to do. ``It makes me look for a proof of this theorem". -- Well, now let's see what your search actually consists in; then we will know, ``how it relates to your belief in the sentence. || what it has to do with believing the sentence."
\end{quote}
(Wittgenstein \emph{Nachlass} Ts-212,XVII-122-19 (WL), (German, 1932--33?)  (IDP). Heinz Wilhelm Kr\"{u}ger, Alois Pichler (2022: inclusion of corrections in XML transcription)) (Accessed 11 Jun 2022)
}
\end{quote}

No doubt the proposal demands a very careful presentation, thorough documentation and a series of justified connections between such a wide spectrum of concepts and techniques, ranging from the philosophical works of Wittgenstein to a fundamental collection of techniques to formalise proofs in mathematics. There has been sustained effort to fulfill a task going back at least 35 years, and the expectation is to continue with the same determination and care as always. No single article has done this and never will, of course. The present paper belongs to a series of articles which began in 1987--1988 (de Queiroz (1987),(1988ab)), followed by some more in the 1990s. This view of proofs and meaning led to a proposal of reformulating intuitionistic type theory\footnote{One of the main points of such a reformulation is concerned with `equality proofs' and the so-called `identity types'. The identity type is arguably one of the most interesting entities of Martin-L\"{o}f type theory (MLTT). From any type $A$, one can construct the identity type $Id_A (x,y)$ whose inhabitants (if any) are proofs of equality between $x$ and $y$. This type establishes the relation of identity between two terms of $A$, i.e., if there is a construction $x =_p y: A$, then $p$ is a witness or proof that $x$ is indeed equal to $y$. Both $Id_A (x,y)$ and $x =_p y: A$ are types of equality, the latter being an explicit way of \emph{definitional} equality, whereas the former is \emph{propositional} equality. The proposal of the Univalence Axiom made the identity type perhaps one of the most studied aspects of type theory in the last decade or so. It proposes that in type theory, to say $x=y$ is equivalent to saying that $x\simeq y$, that is the identity type is equivalent to the type of equivalences. Another important aspect is the fact that it is possible to interpret the paths between two points of the same space. This interpretation gives rise to the view of equality as a collection of homotopical paths. And such connection of type theory and homotopy theory makes type theory a suitable foundation for both computation and mathematics. Nevertheless, in the original formulation of intensional identity type in MLTT, this interpretation is only a semantical one, i.e., all elements are equal to a special term representing reflexivity, and there is no syntactical counterpart for the concept of path in type theory. For that reason, the addition of paths to the syntax of homotopy type theory was proposed by (de Queiroz and Gabbay (1994), de Oliveira and de Queiroz (1999), de Queiroz and de Oliveira (2011,2014),  de Queiroz et al. (2016), Ramos et al. (2017), Ramos et al. (2018), de Veras et al. (2023a)), in these works, the authors use a formal entity called `computational path', proposed by (de Queiroz and Gabbay (1994)), and show that it can be used to formalise the identity type in a more explicit manner. As we have emphasised in (Ramos et al.\ (2021ab), de Veras et al.\ (2023ab)), this allows for making useful bridges between theory of computation, algebraic topology, logic, higher categories, and higher algebra (Mart\'{\i}nez-Rivillas and de Queiroz (2022ab,2023ab)), and a single concept seems to serve as a bridging bond: ``path". This is made possible by using an alternative formulation of the identity type which provides an explicit formal account of ``path", operationally understood as an invertible sequence of rewrites (such as Church's ``conversion" between $\lambda$-terms).}
which gave rise to technical results such as those documented in a book (2011), some PhD theses (e.g., Ramos (2018), Mart\'{\i}nez-Rivillas (2022)), as well as in several articles spanning at least three decades up until 2023. 

\paragraph*{Acknowledgements.} First and foremost, we acknowledge the role of Anderson Nakano (PUC-SP) and Marcos Silva (UFPE), the organisers of the online meeting \emph{100 anos do Tractatus Logico-Philosophicus} (14--17 Sept 2021), who kindly accepted the offer of a talk. This served as additional stimulus to write this paper, which, in turn led us to new findings in Wittgenstein's \emph{Nachlass}. Equally important was the very careful and very positive work done by the anonymous reviewer who has significantly contributed to improve the text.

\section*{References}
Dummett, M. 1977. \textit{Elements of Intuitionism}, Clarendon Press, Oxford.\\
\ \\
Dummett, M. 1991. \textit{The Logical Basis of Metaphysics}, Harvard University Press, Cambridge (Mass.).\\
\ \\
Gabbay, D.M., de Queiroz, R.J.G.B. 1992. ``Extending the Curry-Howard Interpretation to Linear, Relevant and Other Resource Logics" \emph{Journal of Symbolic Logic} 57(4):1319--1365.\\
\ \\
Gentzen, G. 1935. Untersuchungen \"uber das logische Schlie{\ss}en. I. \emph{Math Z} 39, 176--210. https://doi.org/10.1007/BF01201353. (English translation in J.\ von Plato (2017)).\\
\ \\
Henkin, L. 1961. ``Some remarks on infinitely long formulas". In \emph{Infinitistic Methods (Proc. Sympos. Foundations of Math.)}, Warsaw, pages 167--183.
Pergamon, Oxford.\\
\ \\
Hermes, H. 1959. ``Zum Inversionsprinzip der operativen Logik". In Heyting,
A., editor, \emph{Constructivity in Mathematics}, pages 62--68. North-Holland,
Amsterdam.\\
\ \\
Hilmy, S. 1987. \textit{The Later Wittgenstein: The Emergence of a New Philosophical Method}, Blackwell, Oxford.\\
\ \\
Hintikka,  J. 1968. ``Language-games for quantifiers". In N.\ Rescher, editor, \emph{Studies in Logical Theory}, pages 46--72. Blackwell.\\
\ \\
Hintikka,  J. 1979. ``Quantifiers vs. Quantification Theory". In E. Saarinen (ed) \textit{Game-Theoretical Semantics}. Synthese language library, vol 5. Springer, Dordrecht.\\
\ \\
Hintikka, M. and Hintikka, J. 1986. \textit{Investigating Wittgenstein}.  Basil Blackwell, Oxford.\\
\ \\
Hodges, W. 2001. ``A Sceptical Look". \textit{Proceedings of the Aristotelian Society} Supplementary Volume {75}(1):17--32.\\
\ \\
Howard, W. 1980. ``The formulae-as-types notion of construction".
In H. Curry, J.R. Hindley, J. Seldin (eds.), \emph{To H. B. Curry: Essays on Combinatory Logic, Lambda Calculus, and Formalism}. Academic Press (1980). \\
\ \\
Howard, W. 2014. Wadler's Blog, 2014, https://wadler.blogspot.com/2014/08/howard-on-curry-howard.html\\
\ \\
Keiff, L. 2009. ``Dialogical Logic". In \textit{Stanford Encyclopedia of Philosophy}, Stanford. \linebreak https://plato.stanford.edu/entries/logic-dialogical/\\
\ \\
Krabbe, E. 2001. ``Dialogue Logic Restituted". In W.\ Hodges and E.\ C.\ W.\ Krabbe (eds.) \textit{Dialogue Foundations}, Wiley, pp.\ 33--49.\\
\ \\
Lorenzen, P. 1950. ``Konstruktive Begr\"{u}ndung der Mathematik". \emph{Mathematische Zeitschrift} 53(2):162--202.\\
\ \\
Lorenzen, P. 1955. \textit{Einf\"{u}hrung in die operative Logik und Mathematik}. Die Grundlehren der mathematischen Wissenschaften, vol. 78. Springer-Verlag, Berlin-GÃ¶ttingen-Heidelberg. VII + 298 pp\\
\ \\
Lorenzen, P. 1969. \textit{Normative Logic and Ethics\/}, series B.I-Hochschultaschenb\"ucher. Systematische Philosophie, vol. 236$^*$, Bibliographisches Institut, Mannheim/Z\"urich.\\
\ \\
Martin-L\"of, P. 1984. \emph{Intuitionistic Type Theory}. (Notes taken by G.\ Sambin). Bibliopolis, Napoli. ISBN 978-8870881059.\\
\ \\
Martin-L\"of, P. 1987. ``Truth of a Proposition, Evidence of a Judgement, Validity of a Proof". \textit{Synthese} {73}:407--420.\\
\ \\
Martin-L\"of, P.  2013. ``Verificationism Then and Now". Chapter 1 of M.\ van der Schaar (ed.), \emph{Judgement and the Epistemic Foundation of Logic},  Logic, Epistemology, and the Unity of Science 31, DOI 10.1007/978-94-007-5137-8\_1, Springer 2013.\\
\ \\
Martin-L\"of, P.  2019. ``Logic and Ethics". In Piecha, T.; Schroeder-Heister, P. (ed.), \emph{Proof-Theoretic Semantics: Assessment and Future Perspectives}. Proceedings of the Third T\"ubingen Conference on Proof-Theoretic Semantics, 27--30 March 2019,  Univ T\"ubingen, pp.\ 227--235. http://dx.doi.org/10.15496/publikation-35319.\\
\ \\
Mart\'{\i}nez-Rivillas, D.O. 2022. \emph{Towards a homotopy domain theory}. PhD thesis,
CIn-UFPE (November 2022). Centro de Inform\'atica, Universidade Federal
de Pernambuco, Recife, Brazil. https://repositorio.ufpe.br/handle/123456789/49221\\
\ \\
Mart\'{\i}nez-Rivillas, D.O., de Queiroz, R.J.G.B. 2022a. ``$\infty$-Groupoid Generated by an Arbitrary Topological $\lambda$-Model". \emph{Logic J.\ of the IGPL} 30(3):465--488.\\
\ \\
Mart\'{\i}nez-Rivillas, D.O., de Queiroz, R.J.G.B. 2022b. ``Towards a Homotopy Domain Theory". \emph{Archive for Mathematical Logic}. Nov 2022. https://doi.org/10.1007/s00153-022-00856-0\\
\ \\
Mart\'{\i}nez-Rivillas, D.O., de Queiroz, R.J.G.B. 2023a. ``The Theory of an Arbitrary Higher $\lambda$-Model". To appear in \emph{Bulletin of the Section of Logic}. 2023. arXiv:2111.07092\\
\ \\
Mart\'{\i}nez-Rivillas, D.O., de Queiroz, R.J.G.B. 2023b. ``Solving homotopy domain equations". (Submitted for publication.) arXiv:2104.01195\\
\ \\
Moore, G. E. 1955. ``Wittgenstein's Lectures in 1930--33". \textit{Mind} {54}(253):1--27.\\
\ \\
de Oliveira, A.G., de Queiroz, R.J.G.B. 1999. ``A Normalization Procedure for the Equational Fragment of
Labelled Natural Deduction". \emph{Logic J.\ of the IGPL} 7(2):173--215.\\
\ \\
de Oliveira, A.G., de Queiroz, R.J.G.B. 2005. ``A new basic set of proof transformations". In S.\ Artemov, H.\ Barringer, A.\ Garcez, L.\ Lamb, \& J.\ Woods (Eds.), \emph{We will show them! Essays in Honour of Dov Gabbay} (Vol. 2, pp.\ 499--528). London: College Publications.\\
\ \\
Pears, D. 1987. \textit{The False Prison.  A Study of the Development of Wittgenstein's Philosophy}. Volume I. Clarendon Press, Oxford.\\
\ \\
Pears, D. 1990. ``Wittgenstein's Holism". \emph{Dialectica} 44(2):165--173.\\
\ \\
Peirce, C. S. 1932. \emph{Collected Papers of Charles Sanders Peirce, Volumes I and II: Principles of Philosophy and Elements of Logic}. C. Hartshorne and
P.\ Weiss (ed.). Harvard Univ Press.\\
\ \\
Pietarinen, A-V. 2014. ``Logical and Linguistic Games from Peirce to Grice to Hintikka". \textit{Teorema}
{33}(2) 121--136.\\
\ \\
von Plato, J. 2017. \emph{Saved from the Cellar.
Gerhard Gentzen's Shorthand Notes on Logic and Foundations of Mathematics.} Springer.\\
\ \\
Prawitz, D. 1965. \textit{Natural Deduction: A Proof-Theoretical Study}. Acta Universitatis Stockholmiensis, Stockholm Studies in Philosophy no.\ 3. Almqvist \& Wiksell, Stockholm, G\"oteborg, and Uppsala.\\
\ \\
Prawitz, D. 1977. ``Meaning and proofs: on the conflict between classical and intuitionistic logic". \textit{Theoria} (Sweden) {XLIII} 2--40.\\
\ \\
Prawitz, D.  2019. ``Validity of Inferences Reconsidered". In Piecha, T.; Schroeder-Heister, P. (ed.), \emph{Proof-Theoretic Semantics: Assessment and Future Perspectives}. Proceedings of the Third T\"ubingen Conference on Proof-Theoretic Semantics, 27--30 March 2019,  Univ T\"ubingen, pp.\ 213--226. http://dx.doi.org/10.15496/publikation-35319.\\
\ \\
de Queiroz, R.J.G.B. 1987. ``Note on Frege's notions of definition and the relationship proof theory vs.\ recursion theory (Extended Abstract)". In \emph{Abstracts of the VIIIth International Congress of Logic, Methodology and Philosophy of Science.} Vol.\ 5, Part I, Institute of Philosophy of the Academy of Sciences of the USSR, Moscow, 1987, pp.\ 69--73.\\
\ \\
de Queiroz, R.J.G.B. 1988a. ``A Proof-Theoretic Account of Programming and the Role of Reduction Rules". \emph{Dialectica} 42(4):265--282.\\
\ \\
de Queiroz, R.J.G.B. 1988b. ``The mathematical language and its semantics: to show the consequences of a proposition is to give its meaning". In P.\ Weingartner, G.\ Schurz, E.\ Leinfellner, R.\ Haller, A.\ H\"ubner, (eds.) \emph{Reports of the thirteenth international Wittgenstein symposium} 18, pp.\ 259--266.\\
\ \\
de Queiroz, R.J.G.B. 1989. ``Meaning, function, purpose, usefulness, consequences -- interconnected concepts (abstract)". In \emph{Abstracts of Fourteenth International Wittgenstein Symposium (Centenary Celebration)}, 1989, p.\ 20. Symposium held in Kirchberg/ Wechsel, August 13--20 1989.\\
\ \\
de Queiroz, R.J.G.B. 1991. ``Meaning as grammar plus consequences". \emph{Dialectica} 45(1):83--86.\\
\ \\
de Queiroz, R.J.G.B. 1992. ``Grundgesetze alongside Begriffsschrift (abstract)", in \emph{Abstracts of Fifteenth International Wittgenstein Symposium}, 1992, pp.\ 15--16. Symposium held in Kirchberg/Wechsel, August 16--23 1992.\\
\ \\
de Queiroz, R.J.G.B. 1994. ``Normalisation and Language Games". \emph{Dialectica} 48(2):83--123.\\
\ \\
de Queiroz, R.J.G.B. 2001. ``Meaning, function, purpose, usefulness, consequences -- interconnected concepts". \emph{Logic J of the IGPL} 9(5):693--734.\\
\ \\
de Queiroz, R.J.G.B. 2008. ``On Reduction Rules, Meaning-as-use, and Proof-theoretic Semantics". \emph{Studia Logica} 90:211--247.\\
\ \\
de Queiroz, R.J.G.B., Gabbay, D.M. 1994. ``Equality in Labelled Deductive Systems and the functional interpretation of propositional equality". In P.\ Dekker, and M.\ Stokhof, (eds.), \emph{Proceedings of the 9th Amsterdam Colloquium 1994}, ILLC/Department of Philosophy, University of Amsterdam, pp.\ 547--546.\\
\ \\
de Queiroz, R.J.G.B., Gabbay, D.M.  1995. ``The functional interpretation of the existential quantifier"'. \emph{Bull of the IGPL} 3(2--3):243--290.\\
\ \\
de Queiroz, R.J.G.B., Gabbay, D.M.  1997. ``The functional interpretation of modal necessity". In M.\ de Rijke, (ed.), \emph{Advances in Intensional Logic}, Applied Logic Series, Kluwer, 1997, pp.\ 61--91.\\
\ \\
de Queiroz, R.J.G.B., Gabbay, D.M. 1999. Labelled Natural Deduction. In: Ohlbach, H.J., Reyle, U. (eds) \emph{Logic, Language and Reasoning}. Trends in Logic, vol 5. Springer, Dordrecht. https://doi.org/10.1007/978-94-011-4574-9\_10\\
\ \\
de Queiroz, R.J.G.B., Maibaum, T.S.E. 1990. ``Proof Theory and Computer Programming". \emph{Zeitschrift f\"ur mathematische Logik und Grundlagen der Mathematik} 36:389--414. \\
\ \\
de Queiroz, R.J.G.B., Maibaum, T.S.E. 1991. ``Abstract Data Types and Type Theory: Theories as Types". \emph{Zeitschrift f\"ur mathematische Logik und Grundlagen der Mathematik} 37:149--166. \\
\ \\
de Queiroz, R.J.G.B., de Oliveira, A.G. 2011. ``The Functional Interpretation of Direct Computations". \emph{Electronic Notes in Theoretical Computer Science} 269:19--40.\\
\ \\
de Queiroz, R.J.G.B., de Oliveira, A.G. 2014.  ``Natural Deduction for Equality: The Missing Entity". In Pereira, L., Haeusler, E., de Paiva, V. (eds) \emph{Advances in Natural Deduction}. Pages 63--91. Trends in Logic, vol 39. Springer, Dordrecht..\\
\ \\
de Queiroz, R.J.G.B., de Oliveira, A.G., Gabbay, D.M. 2011.  \emph{The Functional Interpretation of Logical Deduction}. Vol. 5 of Advances in Logic series. Imperial College Press / World Scientific, Oct 2011.\\
\ \\
de Queiroz, R.J.G.B., de Oliveira, A.G., Ramos, A.F. 2016. ``Propositional Equality, Identity Types, and Computational Paths". \emph{South Amer. J. of Logic} 2(2):245--296.\\
\ \\
Ramos, A.F. 2018. \emph{Explicit computational paths in type theory}. PhD thesis,
CIn-UFPE (August 2018). Centro de Inform\'atica, Universidade Federal
de Pernambuco, Recife, Brazil. \linebreak https://repositorio.ufpe.br/handle/123456789/32902 (Abstract in: Ramos, A. (2019). Explicit Computational Paths in Type Theory. \emph{Bulletin of Symbolic Logic}, 25(2):213-214. \linebreak doi:10.1017/bsl.2019.2)\\
\ \\
Ramos, A.F., de Queiroz, R.J.G.B., de Oliveira, A.G.  2017. ``On the identity type as the type of computational paths". \emph{Logic J.\ of the IGPL} 25(4):562--584.\\
\ \\
Ramos, A.F., de Queiroz, R.J.G.B., de Oliveira, A.G., de Veras, T.M.L. 2018. ``Explicit Computational Paths". \emph{South Amer.\ J.\ of Logic} 4(2):441--484.\\
\ \\
Ramos, A.F., de Queiroz, R.J.G.B., de Oliveira, A.G. 2021a. ``Computational Paths and the Fundamental Groupoid of a Type". In: \emph{Encontro de Teoria da Computa\c{c}\~ao (ETC)}, 6, Evento Online. Porto Alegre: Sociedade Brasileira de Computa\c{c}\~ao, 2021. p.\ 22--25. ISSN 2595-6116. DOI: https://doi.org/10.5753/etc.2021.16371\\
\ \\
Ramos, A.F., de Queiroz, R.J.G.B., de Oliveira, A.G. 2021b. ``Convers\~ao de Termos, Homotopia, e Estrutura de Grup\'oide". In: \emph{Workshop Brasileiro de L\'ogica (WBL)}, 2, Evento Online. Porto Alegre: Sociedade Brasileira de Computa\c{c}\~ao, 2021. p.\ 33--40. ISSN 2763-8731. DOI: https://doi.org/10.5753/wbl.2021.15776\\
\ \\
Schroeder-Heister, P. 2006. ``Validity concepts in proof-theoretic semantics". \emph{Synthese} 148 (3):525--571.\\
\ \\
Schroeder-Heister, P. 2018. ``Proof-Theoretic Semantics" in \textit{Stanford Encyclopedia of Philosophy}, Stanford University. https://plato.stanford.edu/entries/proof-theoretic-semantics/\\
\ \\
V\"an\"a\"anen, J. 2022. {The Strategic Balance of Games in Logic}. arXiv:2212.01658\\
\ \\
de Veras, T.M.L., Ramos, A.F., de Queiroz, R.J.G.B., de Oliveira, A.G.,  2021. ``Calculation of Fundamental Groups via Computational Paths". In: \emph{Encontro de Teoria da Computa\c{c}\~ao (ETC)}, 6. , 2021, Evento Online. Porto Alegre: Sociedade Brasileira de Computa\c{c}\~ao, 2021 . p.\ 17--21. ISSN 2595-6116. DOI: https://doi.org/10.5753/etc.2021.16370.\\
\ \\
de Veras, T.M.L., Ramos, A.F., de Queiroz, R.J.G.B., de Oliveira, A.G.,  2023a. ``A Topological Application of Labelled Natural Deduction". To appear in \emph{South American J.\ of Logic}. arXiv:1906.09105\\
\ \\
de Veras, T.M.L., Ramos, A.F., de Queiroz, R.J.G.B., Silva, T.D.O., de Oliveira, A.G.,  2023b. ``Computational Paths - A Weak Groupoid". arXiv: 2007.07769. (Submitted for publication)\\
\ \\
Wittgenstein, L. 1974. \textit{Letters to Russell, Keynes and Moore},
Ed. with an Introd. by G.\ H.\ von Wright, (assisted by B.\ F.\ McGuinness), Basil Blackwell, Oxford.\\
\ \\
Wittgenstein, L. 1982. \emph{Last Writings on the Philosophy of Psychology}, vol.\ 1, 1982, vol. 2, 1992, G.H.\ von Wright and H.\ Nyman (eds.), trans. C.G.\ Luckhardt and M.A.E.\ Aue (trans.), Oxford: Blackwell.\\
\ \\
Wittgenstein, L. 1974. \textit{Lectures and Conversations on Aesthetics, Psychology and Religious Belief}, 1966, C.\ Barrett (ed.), Oxford: Blackwell.\\
\ \\
Wittgenstein, L. 1961. \textit{Notebooks 1914--1916,} G.\ H.\ von Wright and G.\ E.\ M.\ Anscombe (eds.), Oxford: Blackwell.\\
\ \\
Wittgenstein, L. 1969. \textit{On Certainty,} G. E. M. Anscombe and G. H. von Wright (eds.), G.E.M.\ Anscombe and D. Paul (trans.), Oxford: Blackwell.\\
\ \\
Wittgenstein, L. 1974. \textit{Philosophical Grammar,} R.\ Rhees (ed.), A.\ Kenny (trans.), Oxford: Blackwell.\\
\ \\
Wittgenstein, L. 1953. \textit{Philosophical Investigations,} G.E.M.\ Anscombe and R.\ Rhees (eds.), G.E.M.\ Anscombe (trans.), Oxford: Blackwell.\\
\ \\
Wittgenstein, L. (1964) \textit{Philosophical Remarks}, R.\ Rhees (ed.), R.\ Hargreaves and R.\ White (trans.), Oxford: Blackwell.\\
\ \\
Wittgenstein, L. 1971. \textit{ProtoTractatus--An Early Version of Tractatus Logico- Philosophicus}, B.\ F.\ McGuinness, T.\ Nyberg, G.\ H.\ von Wright (eds.), D.\ F.\ Pears and B.\ F.\ McGuinness (trans.), Ithaca: Cornell University Press.\\
\ \\
Wittgenstein, L. 1956. \textit{Remarks on the Foundations of Mathematics}, G.H.\ von Wright, R.\ Rhees and G.\ E.\ M.\ Anscombe (eds.), G.\ E.\ M.\ Anscombe (trans.), Oxford: Blackwell, revised edition 1978.\\
\ \\
Wittgenstein, L. 1980. \textit{Remarks on the Philosophy of Psychology,}, vol.\ 1, G. E. M. Anscombe and G. H. von Wright (eds.), G.\ E.\ M.\ Anscombe (trans.), vol. 2, G. H. von Wright and H. Nyman (eds.), C.\ G.\ Luckhardt and M.A.E.\ Aue (trans.), Oxford: Blackwell.\\
\ \\
Wittgenstein, L. 1922. 
\emph{Tractatus Logico-Philosophicus},  Translated by C.K.\ Ogden.
Kegan Paul, Trench, Trubner \& Co., Ltd. New York: Harcourt, Brace \& Company, Inc.\\
\ \\
Wittgenstein, L. 2016. Wittgenstein, Ludwig: Interactive Dynamic Presentation (IDP) of Ludwig Wittgenstein's philosophical \emph{Nachlass} [wittgensteinonline.no]. Edited by the Wittgenstein Archives at the University of Bergen (WAB) under the direction of Alois Pichler. Bergen: Wittgenstein Archives at the University of Bergen 2016-. http://wab.uib.no/wab\_BEE.page (accessed 13 May 2021)

\end{document}